\documentclass[12pt]{article}

\usepackage[matrix,arrow]{xy}
\usepackage{verbatim}
\usepackage{graphicx}

\usepackage[usenames,dvipsnames]{xcolor}

\definecolor{light-gray}{gray}{0.50}

\title{Double Hurwitz numbers and multisingularity loci in genus~0}
\author{Maxim Kazarian\thanks{Steklov Mathematical Institute RAS,
National Research University Higher School of Economics,
Skolkovo Institute of Science and Technology, e-mail: kazarian@mccme.ru},
Sergei Lando\thanks{National Research University Higher School of Economics,
Skolkovo Institute of Science and Technology, e-mail: lando@hse.ru},
Dimitri Zvonkine\thanks{
Universit\'e de Versailles St-Quentin,
CNRS, 45 Avenue des \'Etats Unis, 78000 Versailles, e-mail: dimitri.zvonkine@gmail.fr.
}}

\date{}

\usepackage{amssymb,amsmath}
\usepackage{theorem}

\def\C{{\mathbb C}}
\def\CP{{\mathbb C}{\rm P}}

\def\Q{{\mathbb Q}}
\def\Z{{\mathbb Z}}

\def\D{{\cal D}}

\def\cA{{\cal A}}
\def\cH{{\cal H}}
\def\ocH{{\overline{{\cal H}}}}
\def\cL{{\cal L}}

\def\ocM{{\overline{{\cal M}}}}

\def\cO{{\cal O}}

\def\cR{{\cal R}}

\def\f{{\tilde f}}

\def\deg{{\rm deg}}
\def\Hom{{\rm Hom}}

\def\Aut{{\rm Aut\,}}

\let\a\alpha
\let\b\beta
\let\m\mu
\let\p\partial

\let\s\sigma

\let\k\kappa
\let\L\Lambda
\def\oX{\overline X}
\def\oPsi{\overline\Psi}

\def\hf{\hat{f}}
\def\hz{\hat{z}}
\def\hX{\widehat{X}}
\def\Bl{{\widetilde{P\ocH}_{r+1|\k}}}
\def\opsi{\widehat{\psi}}
\def\oxi{\widehat{\xi}}
\def\tX{\widetilde{\mathcal{X}}}

\def\p{\partial}
\def\h{\hbar}

\newtheorem{theorem}{Theorem}[section]
\newtheorem{lemma}[theorem]{Lemma}

\newtheorem{proposition}[theorem]{Proposition}
\newtheorem{corollary}[theorem]{Corollary}

{\theorembodyfont{\rmfamily}    
\newtheorem{remark}[theorem]{Remark}

\newtheorem{definition}[theorem]{Definition}

}

\def\mdskip{\ifnum\medskipamount>\lastskip
      \vskip-\lastskip\vskip\medskipamount\fi}
\def\QED{\ifmmode\eqno\square\else
{\parfillskip0pt\hfil$\square$\par}\mdskip\fi}

\begin{document}
\def\cX{{\mathcal X}}
\def\ocX{{\overline{\mathcal X}}}
\def\bp{{\bf p}}
\def\bt{{\bf t}}
\let\l\lambda

\maketitle

\begin{abstract}
In the Hurwitz space of rational functions on $\CP^1$ with poles of given orders, we study the loci of multisingularities, that is, the loci of functions with a given ramification profile over~0. We prove a recursion relation on the Poincar\'e dual cohomology classes of these loci and deduce a differential equation on Hurwitz numbers.
\end{abstract}

\tableofcontents

\section{Introduction}

{\em Double Hurwitz numbers} count transitive factorizations of the identity permutation in the symmetric group $S_K$ into a product of two permutations with given cycle types and a given number of transpositions.  In topological terms, they enumerate topologically distinct meromorphic functions $f$ on Riemann surfaces $C$ of a given genus~$g$
with prescribed orders of poles and zeroes and prescribed nonzero simple critical values. In the case when
the genus of $C$ equals $g=0$ and one of the two distinguished permutations is the identity, a closed formula for these numbers was proposed by Hurwitz more than a century ago. In general, in genus~0 explicit formulas for the generating function of double Hurwitz numbers are known. In spite of the existence of such formulas and a variety of modern proofs, many natural questions concerning these numbers remain open.

Hurwitz's argument was algebraic, based on the study of combinatorics of the permutation group.
On the other hand, Hurwitz numbers are related to the geometry of spaces of meromorphic functions and their compactifications, called {\em Hurwitz spaces}.
We propose a new recursion for genus~0 double Hurwitz numbers that has a topological origin: it is derived
from cohomological identities on loci of functions with a given multisingularity, that is, a given ramification profile over one branch point. We expect that variations of
this approach could be adapted to other families of Hurwitz numbers for
which effective formulas are not known at the moment, including those for higher genus curves.

This paper is a continuation of our project on the study of cohomology classes of strata in Hurwitz spaces initiated in~\cite{KL1,KL1a,LZ1,LZ2,KLZ}.

Consider a generic map $f : X \to Y$ between two smooth compact manifolds. Given a singularity type $\alpha$, denote by $X_\alpha$ the locus of points in $X$ where $f$ has a singularity of this type. Thom's principle states that Poincar\'e dual cohomology class of $X_\alpha$ is a universal polynomial in Chern classes of the map $f$ that only depends on $\alpha$.
Further, consider the locus $Y_{\{ \alpha_1, \dots, \alpha_k\}}$ of points in $Y$ whose preimage contains $k$ points of given singularity types. This locus is called a {\em multisingularity locus}. Kazarian's principle states that the Poincar\'e dual cohomology class of a multisingularity locus is the push-forward under $f$ of a universal polynomial in Chern classes of $f$ that only depends on $\alpha_1, \dots, \alpha_k$. The universal map over a Hurwitz space does not satisfy the genericity assumptions for these principles to apply. However it turns out that one can compute the corrections to the universal polynomials that appear in this case.




\subsection{Double Hurwitz numbers}

Given a partition $\l = (\l_1, \l_2, \dots)$ we denote by $\ell(\l)$ the number of its parts, by $|\l|$ the sum of its parts, by $|\Aut \l|$ the number of permutations of its parts that preserve their values, and by $p_\l$ and $q_\l$ the monomials $\prod p_{\l_i}$ and $\prod q_{\l_i}$,
respectively.  For instance, if $\l = (4,3,3,1,1,1)$, then $\ell(\l) =6$, $|\l| = 13$, $|\Aut \l| = 2! \cdot 3!$, $p_\l = p_1^3 p_3^2 p_4$, $q_\l = q_1^3 q_3^2 q_4$.

Let $\l$ and $\mu$ be two partitions of a positive integer~$K$.
Let $g \geq 0$ be an integer and denote $m= \ell(\l) + \ell(\mu) + 2g-2$. Thus $g$ can be recovered from $m$ and vice versa. A {\em transitive factorization} of genus~$g$ and ramification type $(\l, \mu)$ is a list of two permutations $\s, \rho \in S_K$ and $m$ transpositions $\tau_1, \dots, \tau_m \in S_K$ in the symmetric group $S_K$ such that
\begin{itemize}
\item the cycle type of $\s$ is $\l$, the cycle type of $\rho$ is $\mu$;
\item the product $\rho \tau_m \dots \tau_1 \s$ is the identity permutation;
\item the subgroup of $S_K$ generated by $\s$, $\rho$, and the transpositions~$\tau_1,\dots,\tau_m$ is transitive.
\end{itemize}
The {\em double Hurwitz number} $h_{g;\l,\mu}$ is the number of transitive factorizations as above divided by $K!$.

Double Hurwitz numbers can be organized into a natural \emph{generating function} $H$, which is an infinite power series in two infinite sets of variables, $p=(p_1,p_2,\dots)$ and $q=(q_1,q_2,\dots)$, and an additional formal variable~$\b$. The coefficient of the monomial $\b^m p_\lambda q_\mu$ in $H$ is
$$
\frac 1{m! \, | \Aut \lambda| \, |\Aut \mu |}\, h_{g;\lambda, \mu}.
$$

The generating function~$H$ can be determined, for example, by the so-called \emph{cut-and-join equation}~\cite{GJ}
$$\frac{\p}{\p\b}e^H=W\,e^H,\qquad W=\frac12\sum_{i,j}\Bigl(
(i+j)\,p_ip_j\frac{\p}{\p p_{i+j}}+
i\,j\,p_{i+j}\frac{\p^2}{\p p_i\p p_j}
\Bigr).$$
This equation determines the function $H$ uniquely from the initial conditions
 $H(0;p,q)=\sum_{n=1}^\infty\frac{p_nq_n}{n}$. This initial condition indicates that the only degree~$n$ map $C \to \CP^1$ unramified outside~$0$ and~$\infty$ is the map $\CP^1 \to \CP^1$ given by~$z\mapsto z^n$. More explicitly, the series~$e^H$ can be written as
$$e^H=e^{\b\,W}e^{\sum\frac{p_nq_n}{n}}.$$

The Schur polynomials $s_\l(p)$ form an eigenbasis for the operator $W$, and  the corresponding eigenvalues are given by
$$
w(\l) = \frac12\sum_{i=1}^\ell\Bigl(\bigl(\l_i-i+\frac12\bigr)^2-\bigl(-i+\frac12\bigr)^2\Bigr),
$$
where $\l = (\l_1,\l_2,\dots,\l_\ell)$, $\l_1\ge\l_2\ge\dots\ge\l_\ell>0$,
see details, for example, in~\cite{KL4}.
From this one deduces the so-called {\em Frobenius formula} for $H$, see~\cite{LZ,Ok1}:
$$
e^H=\sum_\l e^{w(\l)\,\b}s_\l(p)s_\l(q),
$$
where the summation runs over the set of all Young diagrams (partitions)
$\l$ and $s_\l$ is the corresponding Schur polynomial.

The generating function~$H$ possesses interesting integrable properties. In particular,
it satisfies equations of the Toda lattice hierarchy, see~\cite{Ok1,KL4}.

The \emph{genus expansion} for this function can be obtained by the substitution
$$\h^2H(\h;\h^{-1}\,p;\h^{-1}\,q)=H^{(0)}(p;q)+\h^2 H^{(1)}(p;q)+\dots.$$
Here the term $H^{(g)}$ enumerates genus~$g$ coverings.

Let us collect genus~0 Hurwitz numbers with fixed ramification type~$\l=(\l_1,\dots,\l_r)$ over zero and arbitrary ramification types over infinity into a single series~$h_\l(q)$. More precisely, define the series~$h_\l(q)$ by the expansion
\begin{equation}\label{hdef}
H^{(0)}(p_1+1,p_2,\dots;q)|_{\b=1}=\sum_{r=1}^\infty\frac1{r!}\sum_{\l_1,\dots,\l_r} h_\l(q)\;p_{\l_1}\dots p_{\l_r}.
\end{equation}
Thus the coefficient of $q_\mu$ in $h_\l$ counts the connected genus~0 ramified coverings of the sphere with
\begin{itemize}
\item $r$ marked and numbered zeros of orders $\l_1, \dots, \l_r$;
\item any number of unmarked simple zeros (this is consequence of the shift of variable $p_1$ by 1);
\item a ramification profile $\mu$ over $\infty$;
\item only simple branch points outside $0$ and $\infty$.
\end{itemize}
This coefficient is equal to this number divided by $|\Aut(\mu)|$ and by the factorial of the number of simple branch points.

By definition, the function~$h_\l(q)$ does not depend on the order of the entries $\l_1,\dots,\l_r$.

In this paper we derive a differential equation on the functions $h_\l$ and a broader set of generating functions to be defined later. The equation is of topological origin: it is derived from
cohomological information contained in the stratification of the Hurwitz space by the multisingularity
types of the functions. This differential equation allows one to compute all the double Hurwitz numbers recursively. As far as we can see, this kind of recursion has never appeared before.

Note that at the moment our method meets a serious obstacle in higher genus, due to the fact that, contrary to the genus~0 case,
higher genus Hurwitz spaces are smooth orbifolds no longer.

%
%
%
%
%

\subsection{Hurwitz spaces}

The Hurwitz spaces we consider are spaces of meromorphic functions on a rational curve, with marked poles of prescribed orders. The orders of zeroes are not specified. Functions with given orders of zeroes form certain
subvarieties in the Hurwitz spaces; these subvarieties are the strata of the stratification of these spaces with respect to the multisingularity type.
Since zeros and poles play different roles, the ramification types over~$0$ and~$\infty$ enter our formulas not in a symmetric way.

The derived relations between the singularity strata are independent of the specific Hurwitz space whose stratification we study;
they are totally determined by the local degeneration types of the functions involved. In particular, the $q$-variables, that keep track of the orders of the poles, are not involved explicitly in our equations for the generating function.

For any tuple of positive integers $\k=(k_1,\dots,k_n)$, we denote by
$\cH_\k=\cH_{(k_1,\dots,k_n)}$ the space of rational functions $f(z)$ with~$n$ marked poles of orders $k_1,\dots,k_n$, respectively. Two rational functions are equivalent and correspond to the same point of the Hurwitz space if they can be obtained from each other by a homography $z \mapsto \frac{az+b}{cz+d}$. This Hurwitz space is introduced and discussed in detail in~\cite{ELSV2}. Since we only consider the case when the source curve has genus zero, we do not indicate the genus in the notation.

We'll make use of a generalization of these spaces, namely, for any integer $r\ge1$ we denote by $\cH_{r|\k}=\cH_{r|(k_1,\dots,k_n)}$ the space that parameterizes functions with $n$~marked poles of prescribed orders and~$r$ extra marked points in~$\C P^1$. The marked points that are not poles are called \emph{supplementary marked points}.
We number the marked points in the following order: the supplementary marked points are numbered from $1$ to $r$, while
the poles have numbers $r+1,\dots,r+n$.

We denote by $\ocH_{r|\k}$ a natural completion of the space~$\cH_{r|\k}$. This is the moduli space of stable maps from genus~0 curves to $\CP^1$ with $k$ preimages of $\infty$ of prescribed orders and $r$ more marked points. For more details see~\cite{ELSV2}.

For $r+n \geq 3$, the Hurwitz space is related to the moduli space $\ocM_{r+n}$ of possibly singular stable genus zero curves with~$r+n$ marked points via the natural projection $\pi:\ocH_{r|\k}\to\ocM_{r+n}$ that takes a rational function to the stabilization of its source curve.

This projection is an {\em orbifold cone} (see~\cite{Fulton}, Chapter~4 for the definition of a cone, \cite{ELSV2}, \cite{Sauvaget} for the orbifold generalization).
The group $\C^*$ of nonzero complex numbers acts on $\ocH_{r|\k}$ by multiplying functions by constants. Its projectivization $P\ocH_{r|\k}$, that is, the space of non-constant $\C^*$-orbits, is a compact orbifold of dimension~$\dim P\ocH_{r|\k}=\ell(\k)+|\k|+r-3$. It carries a natural (rational) fundamental homology class
$$[P\ocH_{r|\k}]\in H_{2(\ell(\k)+|\k|+r-3)}(P\ocH_{r|\k})$$
and the characteristic cohomology class of the~$\C^*$-action
$$\xi=c_1(O(1))\in H^2(P\ocH_{r|\k}),$$
where $O(1)$ is the dual of the tautological line bundle of the cone. Here and below we always consider cohomology groups with coefficients in~$\mathbb{Q}$.

\begin{remark} \label{Rem:zerocone}
A particular case that we will encounter in our recursion is $r=n=1$. In this case the moduli space $\ocM_{r+n}=\ocM_2$ does not exist, and we consider $\ocH_{r|k}$ as a cone over the orbifold $\mathcal{B} \Z_k = \{ \mbox{point} \} / \Z_k$. The projectivization $P\ocH_{1|k}$ is a weighted projective space. Similarly, one can view $\ocH_{0|k}$ as a cone over $\mathcal{B} \Z_k$, but we will actually never encounter this space.
\end{remark}

\subsection{Strata and their degrees}

\begin{definition}
For every pure dimension subcone $Z\subset \ocH_{r|\k}$, its \emph{degree}
is defined to be the intersection number of the subvariety~$PZ\subset P\ocH_{r|\k}$ with the suitable power of the class~$\xi$,
$$\deg_{r|\k}(Z)=\int_{PZ}\xi^{\dim PZ}=\int_{P\ocH_{r|\k}}\frac1{1-\xi}\frown[PZ].$$
More generally, for any cohomology class $\a\in H^*(P\ocH_{r|\k})$, we define its degree as the number
$$\deg_{r|\k}(\a)=\int_{P\ocH_{r|\k}}\frac\a{1-\xi}.$$
\end{definition}

The degree is a $\Q$-valued linear function on the total cohomology space of the projectivized Hurwitz space.

In our recursion we will encounter a particular case when $Z$ is included into the zero section of the cone. In this case we set
$$
\deg \, Z =
\begin{cases}
\displaystyle \sum_{z \in Z} \frac1{\mbox{stab}(z)} & \mbox{ if } \dim Z = 0, \\
0 & \mbox{ else.}
\end{cases}
$$

For any tuple of non-negative integers $\l=(\l_1,\dots,\l_r)$, we denote by $\cX_\l=\cX_{(\l_1,\dots,\l_r)}\subset\ocH_{r|\k}$ the set of functions having a zero of order exactly $\l_i$ at the $i$~th supplementary marked point for $i=1,\dots,r$. Thus, if $\l_i=0$, then the function $f$ does not vanish at the $i$~th marked point; if $\l_i=1$ then the function must have a simple zero at the corresponding point, and so on. Further, we denote by $\ocX_\l$ the closure of $\cX_{(\l_1,\dots,\l_r)}$ in $\ocH_{r|\k}$. In particular, with this notation we have $\ocH_{r|\k}=\ocX_{(0,\dots,0)}$, and in general $\ocX_\l$ is a closed suborbifold of $\ocH_{r|\k}$ called a {\em multisingularity stratum}.

Note that we may consider the same type~$\l$ of degeneration of functions in Hurwitz spaces~$\ocH_{r|\k}$ for different~$\k=(k_1,\dots,k_n)$.

\begin{proposition}\label{prophlambda}
For each tuple $\l=(\l_1,\dots,\l_r)$ of positive integers, the generating function for the degrees of the singularity stratum $\ocX_\l$ in different Hurwitz spaces coincides with the generating function~$h_\l(q)$ for the corresponding Hurwitz numbers defined by Eq.~{\rm\eqref{hdef}}:
$$h_\l(q)=\sum_{n=1}^\infty\frac1{n!}\sum_{k_1,\dots,k_n}\deg_{r|(k_1,\dots,k_n)}(\ocX_\l)\;q_{k_1}\dots q_{k_n}.$$
\end{proposition}

This proposition will be proved in Section~\ref{Ssec:HurDeg}. It shows that  computing double Hurwitz numbers is equivalent to computing the degrees of multisingularity strata in Hurwitz spaces.

The Hurwitz space~$\ocH_{r|\k}$ carries natural ``tautological'' classes $\psi_1$, \dots, $\psi_r$. Namely, the cotangent line at the $i$th supplementary marked point to the source curve of a rational function defines a line bundle over $\ocH_{r|\k}$ denoted by~$\cL_i$, and we set $\psi_i=c_1(\cL_i)\in H^2(\ocH_{r|\k})$,
$i=1,\dots,r$. These classes are the pull-backs under the natural projection $\pi:P\ocH_{r|\k}\to\ocM_{r+n}$ of the corresponding $\psi$-classes defined in a similar way on~$\ocM_{r+n}$.

To each singularity stratum $\ocX_\l$ one can assign the fundamental class $[P\ocX_\l]$ of its projectivization and its ``descendants'', obtained by multiplying $[P\ocX_\l]$ by monomials in the $\psi$-classes $\psi_1, \dots, \psi_r$.

\begin{definition}
For any two tuples of non-negative integers $\l=(\l_1,\dots,\l_r)$ and $\nu=(\nu_1,\dots,\nu_r)$,
of the same length $\ell(\l)=\ell(\nu)=r$, we define the function
$$x_{\l,\nu}(q)=\sum_{n=1}^\infty\frac1{n!}\sum_{\k=(k_1,\dots,k_n)}
\deg_{r|\k}(\psi_1^{\nu_1}\dots\psi_r^{\nu_r}\,[P\ocX_\l])\;q_{k_1}\dots q_{k_n}.$$
\end{definition}

The function $x_{\l,\nu}(q)$ is invariant under simultaneous permutations of the parts in~$\l$ and~$\nu$. We have
$$
h_\l(q)=x_{\l,(0,\dots,0)}(q).
$$

\subsection{The recursion}

To formulate the recursion, let us introduce a new set of variables $t_{i,j}$, $i\ge0$, $j\ge0$, and collect the functions $x_{\l,\nu}(q)$ into the following \emph{descendant genus zero Hurwitz potential}
$$X(t,q)=\sum_{r=1}^\infty\frac1{r!}
\sum_{\substack{\l=(\l_1,\dots,\l_r)\\\nu=(\nu_1,\dots,\nu_r)}}
x_{\l,\nu}(q)\;t_{\l_1,\nu_1}\dots t_{\l_r,\nu_r}.$$
The original generating function for the genus zero double Hurwitz numbers is related to the series $X(t,q)$ as follows:
$$H^{(0)}(p_1+1,p_2,\dots;q)|_{\b=1}=X(t,q)\bigm|_{\substack{t_{i,0}=p_i~(i\ge1),\\t_{i,j}=0~(j>0~or~i=0)}}.$$

\begin{proposition} \label{Prop:Init}
We have
$$x_{(0,\dots,0),(\nu_1,\dots,\nu_r)}(q)=\binom{|\nu|}{\nu_1,\dots,\nu_r}z_{|\nu|,r}(q)$$
where
$$z_{d,r}(q)=\sum_{K,n}\frac1{n!}\binom{n+r-3}{d}K^{n+r-3-d}
\sum_{k_1+\dots+k_n=K}
 \prod_{i=1}^n\frac{k_i^{k_i}}{k_i!}q_{k_i}.$$
\end{proposition}

In~\cite{Z1}, an algebra of power series in a single variable has been introduced.
This algebra is generated by two power series, namely,
$$
\sum_{n=1}^\infty n^{n-1}\frac{x^n}{n!},\qquad \sum_{n=1}^\infty n^{n}\frac{x^n}{n!}.
$$
This algebra contains all the generating functions for Hurwitz numbers
enumerating ramified coverings of the sphere with fixed degenerate ramification,
provided the degree~$n$ of the coverings varies. For ramified coverings
of the torus, a similar role is played by the algebra of quasimodular forms.

The power series $z_{d,r}(q)$ play a similar role in
computing double Hurwitz numbers. Namely, all the series $h_\l(q_1,q_2,\dots)$ and $x_{\l,\m}(q_1,q_2,\dots)$
belong to the algebra generated by the series $z_{d,r}$.
Note that the series~$z_{d,r}$ are not algebraically independent, they obey certain polynomial relations.
We will give a description of this algebra in a separate paper.

The expression in Proposition~\ref{Prop:Init} gives the initial conditions for our recursion.
To formulate the recursion itself, we need the following explicit series in the $t$-variables:
$$\Psi_{a,\ell}=\sum_{k=0}^\infty\frac1{k!}
 \sum_{\substack{\nu_1+\dots+\nu_k=\ell+k-3\\\l_1+\dots+\l_k=a}}
 \binom{|\nu|}{\nu_1,\dots,\nu_k}\prod_{i=1}^kt_{\l_i,\nu_i}.$$

 \begin{remark} \label{Rem:Psi}
 The coefficients of $\Psi$ are equal to
 $$
 \binom{|\nu|}{\nu_1,\dots,\nu_k}
 = \int_{\ocM_{\ell+k}} \psi_1^{\nu_1} \cdots \psi_k^{\nu_k}.
 $$
 This is the form in which this series appears in the proof of the recursion.
 \end{remark}

\begin{theorem} \label{Thm:Rec}
The series $X$ obeys the following differential equations valid for any $s\ge0$ and $m\ge0$:
$$\frac{\p X}{\p t_{s+1,m}}=\frac{\p X}{\p t_{s,m}}+s\,\frac{\p X}{\p t_{s,m+1}}-\sum_{\ell\ge1,a\ge0}\frac1{\ell!}\frac{\p\Psi_{a,\ell}}{\p t_{s,m}}
\sum_{\s_1+\dots+\s_\ell=a}\prod_{i=1}^\ell\Bigl(\s_i\,\frac{\p X}{\p t_{\s_i,0}}\Bigr).$$
\end{theorem}

These differential equations provide a recursion for the coefficients of the series~$X$. The recursion expresses each series $x_{\l,\nu}(q)$ with $\ell(\l)=\ell(\nu)=\ell$ as a polynomial in the known series $z_{d,r}(q)$ of Proposition~\ref{Prop:Init}~with $r\le\ell$.

The geometric meaning of the recursion will be explained in Section~\ref{cohomrel}.

\subsection{Examples}

The recursion of Theorem~\ref{Thm:Rec}
translates into \emph{polynomial} relations between the Taylor coefficients $x_{\l,\nu}$ of the $t$-expansion of~$X$. Each series $x_{\l,\nu}$ is expressed as a polynomial in similar series of smaller degrees with respect to the grading defined by $\deg\, x_{\l,\nu}=|\l|$. Thus, due to the initial conditions of the recursion, the series $x_{\l,\nu}$ are expressed as polynomials in the series $z_{d,r}$ of Proposition~\ref{Prop:Init}. Here are a few first examples for the Hurwitz generating series $h_\l(q)=x_{\l,(0,\dots,0)}(q)$:
\begin{align*}
h_{(2)}&=z_{0,1}+z_{1,1},\\
h_{(3)}&=z_{0,1}-\frac{z_{0,1}^2}{2}+3\,z_{1,1}+2\,z_{2,1},\\
h_{(4)}&=z_{0,1}-\frac{5\,z_{0,1}^2}{2}+6\,z_{1,1}-2\,z_{0,1}\,z_{1,1}+11\,z_{2,1}+6\,z_{3,1},\\
h_{(5)}&=z_{0,1}-\frac{15\,z_{0,1}^2}{2}+\frac{5\,z_{0,1}^3}{6}+10\,z_{1,1}-15\,z_{0,1}\,z_{1,1}
 -2\,z_{1,1}^2+35\,z_{2,1}-6\,z_{0,1}\,z_{2,1}\\&\qquad\qquad+50\,z_{3,1}+24\,z_{4,1},\\
h_{(2,2)}&=-6\,z_{0,1}+z_{0,1}^2+z_{0,2}-11\,z_{1,1}+2\,z_{1,2}-6\,z_{2,1}+2\,z_{2,2},\\
h_{(3,2)}&=-10\,z_{0,1}+9\,z_{0,1}^2+z_{0,2}-z_{0,1}\,z_{0,2}-35\,z_{1,1}+6\,z_{0,1}\,z_{1,1}
 +4\,z_{1,2}-z_{0,1}\,z_{1,2}\\&\qquad\qquad-50\,z_{2,1}+8\,z_{2,2}-24\,z_{3,1}+6\,z_{3,2},\\
h_{(2,2,2)}&=85\,z_{0,1}-40\,z_{0,1}^2-18\,z_{0,2}+6\,z_{0,1}\,z_{0,2}+z_{0,3}+225\,z_{1,1}
 -24\,z_{0,1}\,z_{1,1}-51\,z_{1,2}\\&\qquad+6\,z_{0,1}\,z_{1,2}+3\,z_{1,3}+274\,z_{2,1}
 -84\,z_{2,2}+6\,z_{2,3}+120\,z_{3,1}-54\,z_{3,2}+6\,z_{3,3}.
\end{align*}


\subsection{Acknowledgements}

The first two authors are grateful to the participants of the seminar
``Characteristic classes and intersection theory'' at the Department of mathematics,
NRU HSE.
The third author is partly supported by the ANR-18-CE40-0009 ENUMGEOM grant.

\section{Degrees of strata in Hurwitz Spaces} \label{Sec:degrees}

\subsection{Degrees of strata are double Hurwitz numbers} \label{Ssec:HurDeg}

In this section we prove Proposition~\ref{prophlambda}, i.e., the equality
$$
h_\l(q)=\sum_{n=1}^\infty\frac1{n!}\sum_{k_1,\dots,k_n}\deg_{r|(k_1,\dots,k_n)}(\ocX_\l)\;q_{k_1}\dots q_{k_n}.
$$

\emph{Proof}. The argument is similar to that of~\cite{ELSV2}. Let~$m$ be the number of simple critical points of a generic rational function in the stratum~$X_\l$, those where the critical values differ from zero and infinity. This number can be determined by the Riemann--Hurwitz formula. Then there is a natural map from $X_\l$ to $\C^m$ that associates to a rational function $f \in X_\l$ the unordered collection of its~$m$ critical values. This map extends to an algebraic map $\ocX_\l\to\C^m$ called the \emph{Lyashko--Looijenga map}, or \emph{branching morphism}. The target space~$\C^m$ is the $m$th symmetric power of the complex line~$\C$. By definition, the Hurwitz number is the degree of this map, that is, the number of preimages of a generic point. The Lyashko--Looijenga map is, of course, $\C^*$-equivariant, where $\C^*$ acts on the target $\C^m$ with weights $1, \dots, m$. It follows that the Lyashko-Looijenga map descends to a map of projectivizations $\L:P\ocX_\l \to P \C^m$, where $P \C^m$ is the weighted projective space with weights $1, \dots, m$. The map $\L$ has the same degree as the Lyashko-Looijenga map itself.
Moreover, we have $\L^* O(1) = O(1)$ and $\L^* \xi = \xi$. Thus the degree of the Lyashko-Looijenga map or of $\L$ is equal to
$$\deg(\L)=\frac{\int_{P\ocX_\l}\xi^{m-1}}{\int_{P\C^m}\xi^{m-1}}=m!\;\deg_{\ell(\l)|\k}(X_\l).$$
The last equality is due to the fact that  $\int_{P\C^m}\xi^{m-1}=\frac1{m!}$.
We conclude that up to the factor~$m!$, the Hurwitz number (the left-hand side of the equality) is equal to the degree of the singularity stratum~$X_\l$. The factorial factor is accounted for automatically when the numbers are collected into the appropriate generating series. As a result, the two generating series for the Hurwitz numbers and for the degrees of the singularity strata coincide.

In the exceptional case $n=1$, $r=1$, $k_1=\l_1=k$ the stratum $\ocX_\lambda$ is composed of only one map $f = x^k$ and we have $m=0$. The corresponding Hurwitz number is equal to $1/k$, which is also the degree of the cone according to Remark~\ref{Rem:zerocone}.
\QED

\subsection{The initial conditions of the recursion}\label{Ssec:desc}

In this section we compute the degree of the cohomology class $\psi_1^{\nu_1}\dots\psi_r^{\nu_r}[P\ocX_{(0,\dots,0)}]$. Using this, we determine the generating function $x_{(0,\dots,0),\nu}(q)$ and thus prove Proposition~\ref{Prop:Init}.

Recall that $P\ocX_{(0,\dots,0)} = P\ocH_{r|\k}$.
Applying the projection formula to the forgetful map $\pi:P\ocH_{r|\k} \to \ocM_{r+n}$ we get
$$\deg_{r|\k}(\psi_1^{\nu_1}\dots\psi_r^{\nu_r}[P\ocX_{(0,\dots,0)}])=
 \int_{P\ocH_{r|\k}}\frac{\psi_1^{\nu_1}\dots\psi_r^{\nu_r}}{1-\xi}=
 \int_{\ocM_{r+n}}\psi_1^{\nu_1}\dots\psi_r^{\nu_r}\pi_*\Bigl(\frac1{1-\xi}\Bigr).$$
The push-forward class $\pi_*\bigl(\frac1{1-\xi}\bigr)$ is the so-called \emph{total Segre class} of the cone $\ocH_{r|\k}\to\ocM_{r+n}$. It is computed in~\cite{ELSV2}:
$$\pi_*\Bigl(\frac1{1-\xi}\Bigr)=
 \frac{\prod_{i=1}^n\frac{k_i^{k_i}}{k_i!}}{(1-k_1\psi_{r+1})\dots(1-k_n\psi_{r+n})}.$$
Therefore, in this case, the degree of the cohomology class is expressed as an intersection number of $\psi$-classes on~$\ocM_{r+n}$:
\begin{align*}
\int_{\ocM_{r+n}}\psi_1^{\nu_1}\dots\psi_r^{\nu_r}\pi_*\bigl(\frac1{1-\xi}\bigr)&=
\prod_{i=1}^n\frac{k_i^{k_i}}{k_i!}\int_{\ocM_{r+n}}\frac{\psi_1^{\nu_1}\dots\psi_r^{\nu_r}}
{(1-k_1\psi_{r+1})\dots(1-k_n\psi_{r+n})}\\&=
\prod_{i=1}^n\frac{k_i^{k_i}}{k_i!}\binom{n+r-3}{\nu_1,\dots,\nu_r}(k_1+\dots+k_n)^{n+r-3-|\nu|}.
\end{align*}
The last equality is due to the well-known formula for the intersection numbers of $\psi$-classes over $\ocM_{n}$~(see \cite{W1}):
$$\int_{\ocM_{n}}\psi_1^{k_1}\dots\psi_n^{k_n}=
\begin{cases}\binom{n-3}{k_1,\dots,k_n} & \mbox{ if } k_1+\dots+k_n=n-3,\\
0 & \mbox{ if } k_1+\dots+k_n\ne n-3.\end{cases}$$

It follows that
$$
x_{0,\nu}(q)=\sum_{n=1}^\infty\frac1{n!}\sum_{\k=(k_1,\dots,k_n)}
\deg_{r|\k}(\psi_1^{\nu_1}\dots\psi_r^{\nu_r}\,[P\ocX_\l])\;q_{k_1}\dots q_{k_n}
$$
$$
= \binom{|\nu|}{\nu_1, \dots, \nu_r} \sum_{K,n}
\frac1{n!} \binom{n+r-3}{|\nu|} K^{n+r-3-|\nu|}
\sum_{k_1 + \dots + k_n = K} \prod_{i=1}^n
\frac{k_i^{k_i}}{k_i!} q_{k_i}.
$$
This is exactly the formula announced in Proposition~\ref{Prop:Init}. \QED

\section{From a cohomological relation to the recursion}\label{cohomrel}

As explained in the introduction, our recursion formula is a consequence of an identity in the cohomology ring of the Hurwitz space. In this section we introduce and prove this identity, then deduce the recursion formula.



\subsection{The cohomological identity}

Consider the multisingularity variety $P\ocX_{\l_1,\l_2,\dots,\l_r}\subset P\ocH_{r|\k}$. Meromorphic functions forming this variety vanish at the first marked point together with their first $\l_1-1$ derivatives. The subvariety $P\ocX_{\l_1+1,\l_2,\dots,\l_r}\subset P\ocX_{\l_1,\l_2,\dots,\l_r}$ is distinguished by the condition that the next derivative, of order~$\l_1$, also vanishes. We will now show that the derivative of order $\l_i$ is actually a section of some specific line bundle and its vanishing locus can be related to the first Chern class of this line bundle.

Let $C$ be a curve, $f:C\to\C P^1$ a meromorphic function, and $x_1\in C$  a marked point that is not a pole of~$f$. The quotient of the space of $s$-jets at $x_1$ by the space of $(s-1)$-jets at $x_1$ is the line $T_{x_1}^{\otimes s}C$. Thus, if
$f(x_1) = \dots = f^{s-1}(x_1)=0$, then $f^{(s)}(x_1)$
is a linear map $f^{(s)}(x_1):T_{x_1}^{\otimes s}C\to T_0\C P^1$. Note that this map is well defined only if all the derivatives of order less than~$s$ vanish at~$x_1$.

As the function~$f$ varies in the Hurwitz space, the tangent lines $T_{x_1}C$
to the source curves at~$x_1$ form the line bundle $L_1^\vee$ whose first Chern class is $-\psi_1$. Similarly, for the \emph{projectivized} Hurwitz space, the tangent lines $T_0\C P^1$
to the target curve at~$0$ form the line bundle $O(-1)$ with the first Chern class~$-\xi$.
As a corollary, we obtain the following statement.

\begin{lemma} Consider the line bundle $\Hom(\cL_1^{\vee\otimes \l_1},\cO(1))$ over $P\cH_{r|\k}$,
whose first Chern class is $c_1(\Hom(\cL_1^{\vee\otimes \l_1},\cO(-1)))=\l_1\psi_1-\xi$.
Then $f^{(\l_1)}(x_1)$ {\rm(}the derivative of the rational function of order $\l_1$ at the first marked point{\rm)} is a well-defined holomorphic section of this line bundle over the multisingularity stratum $P\ocX_{\l_1,\l_2,\dots,\l_r}$.
\end{lemma}

The zeroes of the holomorphic section $f^{(\l_1)}(x_1)$ form a divisor in $P\ocX_{\l_1,\l_2,\dots,\l_r}$. The multisingularity stratum $P\ocX_{\l_1+1,\l_2,\dots,\l_r}$ is one of the components of this divisor. We will see that the section $f^{(\l_1)}(x_1)$ has a simple vanishing along this component. In general, however, the zero divisor has many other components formed by functions with singular source curves. We call these components the \emph{boundary components} of the zero divisor. Denoting these components by $D_\a$ and the vanishing order of the section along $D_\a$ by $m_\a$, we obtain a cohomological relation
\begin{equation}\label{Dexp}
(\l_1\psi_1 - \xi)[P\ocX_{\l_1,\l_2,\dots,\l_r}]=[P\ocX_{\l_1+1,\l_2,\dots,\l_r}]+\sum m_\a D_\a.
\end{equation}
Our goal is now to list all the components $D_\a$ of the zero divisor and determine the vanishing orders~$m_\a$.

The irreducible components $D_\a$ of the zero divisor correspond to certain singularity types of the functions in the
family $\ocH_{r|\k}$. We start by recalling the description of singularity types from~\cite{KLZ} adapted to our situation,
the one in the presence of supplementary marked points.



\begin{definition} Let $J\subset \{1,\dots,r\}$ and $\ell$ be a positive integer such that $|J|+\ell \geq 3$. Let $\s_1, \dots, \s_\ell$ be positive integers. We denote by $\cX_{\l_1,\l_2,\dots,\l_r}^{J; \s_1, \dots, \s_\ell}$
the locus of meromorphic functions $f\in \ocH_{r|\k}$ satisfying the following conditions. (i)~The source curve~$C$ has an irreducible component~$C_0$, called the {\em central component}, containing the supplementary marked points of the set $J$ and no other marked points; the function $f$ vanishes identically on~$C_0$. (iii)~The component $C_0$ meets the other irreducible components at $\ell$ points; the function~$f$ has zeros of orders exactly $\s_1, \dots, \s_\ell$ at these points. (iii)~The function~$f$ has zeros of orders exactly $\l_i$ at the supplementary marked point that are not in~$J$. We denote by $\ocX_{\l_1,\l_2,\dots,\l_r}^{J; \s_1, \dots, \s_\ell}$ the closure of this locus.
\end{definition}

The closure $\ocX_{\l_1,\l_2,\dots,\l_r}^{J; \s_1, \dots, \s_\ell}$ is a suborbifold of $P\ocH_{r|\k}$ of pure codimension $1+\sum \s_i$. It has many irreducible components corresponding to all possible distributions of the poles and the remaining supplementary marked points among the $\ell$ components of the curve different from the central one. 
If a function~$f$ lies in $\ocX_{\l_1,\l_2,\dots,\l_r}^{J; \s_1, \dots, \s_\ell}$, then its derivatives of \emph{all} orders vanish on the central component. Therefore, functions with nonisolated singularities can appear as the components~$D_\a$ of~\eqref{Dexp} if $1\in J$. The precise statement is given in the following theorem. Denote $|\s|=\sum_{i=1}^\ell \s_i$, $|\l_J| = \sum_{i \in J} \l_i$.

\begin{theorem}\label{cohth} \ \\
{\bf a)} The stratum $P\ocX_{\l_1,\l_2,\dots,\l_r}^{J; \s_1, \dots, \s_\ell}$ lies in the closure of $P\ocX_{\l_1,\dots,\l_r}$ iff $|\s| \geq|\l_J|$. \\
{\bf b)} The stratum $P\ocX_{\l_1,\l_2,\dots,\l_r}^{J; \s_1, \dots, \s_\ell}$
has codimension~$1$ in~$P\ocX_{\l_1,\dots,\l_r}$ iff $|\s|=|\l_J|$. \\
{\bf c)} Assuming this condition is satisfied, the divisor
$P\ocX_{\l_1,\l_2,\dots,\l_r}^{J; \s_1, \dots, \s_\ell}$
is a component of the zero locus of $f^{(\l_1)}(x_1)$ iff $1 \in J$. \\
{\bf d)} In this case the vanishing order of $f^{(\l_1)}(x_1)$
along  $P\ocX_{\l_1,\l_2,\dots,\l_r}^{J; \s_1, \dots, \s_\ell}$ equals $\s_1 \cdots \s_\ell$. \\
{\bf e)}
The vanishing order of $f^{(\l_1)}(x_1)$ along $P\ocX_{\l_1+1,\l_2\dots,\l_r}$ equals~$1$. \\
{\bf f)} The section $f^{(\l_1)}(x_1)$ has no other zeros.
\end{theorem}

\begin{corollary} \label{Cor:coh}
The following relation holds in the cohomology ring of $P\ocH_{r|\k}$:
\begin{multline}\label{Irel}
[P\ocX_{\l_1+1,\l_2,\dots,\l_r}]=(\l_1\psi_1-\xi)[P\ocX_{\l_1,\l_2,\dots,\l_r}]\\-
   \sum_{J,\;1\in J\subset\{1,\dots,r\}}\sum_{\s\vdash|\l_J|}
 \s_1\dots\s_{\ell(\s)}\;[P\ocX_{\l_1,\l_2,\dots,\l_r}^{J; \s_1, \dots, \s_{\ell(\s)}}].
\end{multline}
\end{corollary}

This relation is the precise version of~\eqref{Dexp}.

\subsection{Vanishing loci of $f^{\l_1}(x_1)$: proof of Theorem~\ref{cohth}}

The statement and the proof are similar to~\cite{KLZ}, Section~2.1, but the presence of supplementary points introduces some extra complications. To prove the theorem we first need to discuss deformations of stable maps in order to give a precise description of the closure $P\ocX_{\l_1,\l_2,\dots,\l_r}$ in $P\ocH_{r|\k}$.

Let $f \in \ocH_{r|\k}$ be a stable map. For any point $y \in \C$, a connected component $X$ of $f^{-1}(y)$ that is not a regular point of~$f$ is called a {\em ramification locus}. Thus a ramification locus is either a ramification point of~$f$, or a node of the source curve~$C$ of~$f$, or a stable genus~$0$ subcurve of~$C$ contracted to~$y$. In~\cite{KLZ} we introduced standard deformation spaces for these three types of ramification loci, parametrizing all functions with the same monodromy as the ramification locus. Deformations of ramification points are given by polynomials, deformations of nodes by Laurent polynomials, and deformations of contracted curves by rational functions on contracted curves, respectively. We will recall more precise definitions later.

\begin{proposition} \label{Prop:deformations}
The neighborhood of $f$ in $\ocH_{r|\k}$ is canonically decomposed into a direct product of deformation spaces of the ramification loci of~$f$.
\end{proposition}

\paragraph{Proof.} Draw a disk $D_y$ around every branch point $y \in \C$ of~$f$ so that the disks do not intersect. A sufficiently small deformation $\f$ of~$f$ has all of its branch points in the union of the disks~$D_y$, while over $\C \setminus \cup D_y$ the maps $f$ and $\f$ determine the same nonramified covering. Each connected component of $\f^{-1}(D_y)$, except disks that map biholomorphically to $D_y$, is a small deformation of the corresponding ramification locus~$X \subset f^{-1}(y)$.

Conversely, choose a deformation of every ramification locus of~$f$ sufficiently small for all branch points to lie inside~$D_y$. Then the monodromies around $D_y$ of a ramification locus of~$f$ and of its deformation coincide. Thus the deformations can be glued into the graph of~$f$ instead of the preimage $f^{-1}(D_y)$, and we obtain a stable map in the neighborhood of~$f$. \QED

Let $X$ be a ramification locus of~$f$ over~$0$, that is, $X \subset f^{-1}(0)$. Denote by $\l(X)$ the sum of integers $\l_i$ over the supplementary marked points~$i$ contained in~$X$. Further, denote by $\s(X)$ the total degree of $f$ in the neighborhood of~$X$. More precisely, if $X$ is an isolated zero of order $\l$, then $\s(X)= \l$. If $X$ is a node and $f$ has ramification orders $\s_1$, $\s_2$ at the branches, then $\s(X) = \s_1+\s_2$. If $X$ is a union of irreducible components of~$C$ that meet other components at $\ell$ points at which $f$ has zeros of orders $\s_1, \dots, \s_\ell$, then $\s(X)  = \sum_{j=1}^\ell \s_j$.

\begin{definition}
A ramification locus $X \subset f^{-1}(0)$ is {\em acceptable} if $\lambda(X) \leq \sigma(X)$.
\end{definition}

\begin{lemma} \label{Lem:closure}
A stable map $f$ lies in the closed stratum $P\ocX_{\l_1,\l_2,\dots,\l_r}$ iff (i)~we have $f(x_i)=0$ for all $i$ such that $\lambda_i \geq 1$, and (ii)~all ramification loci in $f^{-1}(0)$ are acceptable.
\end{lemma}

\paragraph{Proof.} First suppose that $f$ lies in  $P\ocX_{\l_1,\l_2,\dots,\l_r}$. Take a smooth map $\hf \in P\cX_{\l_1,\l_2,\dots,\l_r}$ close enough to $f$. The connected components of $\hf^{-1}(D_y)$ are in a one to one correspondence with those of $f^{-1}(y)$ for any branch point~$y$. Moreover, if $X$ is a connected component of $f^{-1}(0)$ and $\hX$ the corresponding connected component of $\hf^{-1}(D_0)$, then $\hX$ contains the same set $J \subset \{1, \dots, n \}$ of supplementary marked points as $X$ and the degree of $\hf$ on $\hX$ is equal to $\s(X)$. It follows that $\s(X) \geq \l(X)$. Indeed, the sum of orders of zeros of $\hf$ over $\hX$ cannot exceed its degree. This proves the ``only if'' part.

To prove the ``if'' part, we will use Proposition~\ref{Prop:deformations}.
First suppose that every ramification locus $X$ of $f$ is either a point or an irreducible curve. If $X \subset f^{-1}(0)$ is a curve, it carries a function $g_X$ with zeros of orders $\lambda_i$ at the supplementary marked points contained in~$X$ (other zeros are allowed) and with poles of orders $\s_j$ at the intersection points of $X$ with other components of the curve~$C$ (no other poles are allowed). The acceptability condition ensures that such a function exists. The functions $g_X$ provide a deformation of~$X$. If $X \in f^{-1}(0)$ is a supplementary marked point, it can be deformed by a polynomial with a root of order exactly $\l_i$ at the marked point. For all other ramification loci (those that correspond to branch points $y \not= 0$ and those that do not contain supplementary marked points) choose any deformation. By Proposition~\ref{Prop:deformations}, this family of deformations determines a deformation of $f$ into the stratum $P\cX_{\l_1,\l_2,\dots,\l_r}$.

Finally, if a ramification locus $X$ is not irreducible, we can consider a family of smooth curves $X_s$ that tend to $X$ as $s \to 0$. The stable map obtained from $f$ by replacing $X$ with $X_s$ lies in the closure $P\ocX_{\l_1,\l_2,\dots,\l_r}$ for any nonzero~$s$, hence $f$ lies in the closure, too.
\QED

\paragraph{Proof of Theorem~\ref{cohth}.}
 \ \\
{\bf a)} {\em The stratum $P\ocX_{\l_1,\l_2,\dots,\l_r}^{J; \s_1, \dots, \s_\ell}$ lies in the closure of $P\ocX_{\l_1,\dots,\l_r}$ iff $|\s| \geq|\l_J|$.} \\
Consider a stable map~$f \in P\cX_{\l_1,\l_2,\dots,\l_r}^{J; \s_1, \dots, \s_\ell}$. It follows immediately from the definition of $P\cX_{\l_1,\l_2,\dots,\l_r}^{J; \s_1, \dots, \s_\ell}$ that the source curve~$C$ of $f$ has only one contracted component, namely the central component of the stratum. This component is acceptable iff $|\s| \geq|\l_J|$. Thus (a) follows from Lemma~\ref{Lem:closure}.

\bigskip

\noindent
{\bf b)} {\em The stratum $P\ocX_{\l_1,\l_2,\dots,\l_r}^{J; \s_1, \dots, \s_\ell}$
has codimension~$1$ in~$P\ocX_{\l_1,\dots,\l_r}$ iff $|\s|=|\l_J|$.}\\
Once we know that $P\ocX_{\l_1,\l_2,\dots,\l_r}^{J; \s_1, \dots, \s_\ell}$ lies in $P\ocX_{\l_1,\dots,\l_r}$, a simple dimension count provides the codimension.

\bigskip

\noindent
{\bf c)} {\em Assuming this condition is satisfied, the divisor
$P\ocX_{\l_1,\l_2,\dots,\l_r}^{J; \s_1, \dots, \s_\ell}$
is a component of the zero locus of $f^{(\l_1)}(x_1)$ iff $1 \in J$.} \\
This is because $f$ has a zero of order exactly $\l_1$ at~$x_i$ when $i \not\in J$, while it vanishes identically in the neighborhood of $x_i$ when $i \in J$.

\bigskip

\noindent
{\bf d)} {\em In this case the vanishing order of $f^{(\l_1)}(x_1)$
along  $P\ocX_{\l_1,\l_2,\dots,\l_r}^{J; \s_1, \dots, \s_\ell}$ equals $\s_1 \cdots \s_\ell$.} \\
The proof repeats that of Lemma~2.1 in~\cite{KLZ} with a slight modification due to supplementary marked points. We give here a rather concise presentation, referring to~\cite{KLZ} for more details.

Let $f$ be a generic stable map in $P\ocX_{\l_1,\l_2,\dots,\l_r}^{J; \s_1, \dots, \s_\ell}$. In particular, the central component $C_0$ of the source curve of $f$ is smooth, the supplementary marked points are zeros of orders exactly $\l_i$, the other zeros are simple, and all branch points of~$f$ except $0$ are simple. According to Proposition~\ref{Prop:deformations}, there is a local chart on the Hurwitz space $\ocH_{r|\k}$ in the neighborhood of $f$ that is a product of standard deformation families for the ramification loci of~$f$. We are going to determine a parametrization of $P\ocX_{\l_1,\l_2,\dots,\l_r}$ in this chart. By the genericity assumption, any ramification locus over $y \not= 0$ is just a simple ramification point. Its deformation family has one parameter corresponding to moving $y$ in~$\C$. Any such deformation keeps $f$ in~$P\ocX_{\l_1,\l_2,\dots,\l_r}$. Every supplementary marked point $x_i$ with $\l_i \geq 2$ and $i \not\in J$ is a ramification locus. The corresponding deformation family is the space of polynomials of degree $\l_i$. Any such deformation is transversal to the stratum $P\ocX_{\l_1,\l_2,\dots,\l_r}$, since it does not preserve the zero of order $\l_i$ at~$x_i$. Thus, in the parametrization of the stratum $P\ocX_{\l_1,\l_2,\dots,\l_r}$, all the parameters of the deformations must be set to~$0$. Finally, the most important ramification locus is the central component $X = C_0$. Its deformation family has {\em moduli parameters}, corresponding to local coordinates in the moduli space $\ocM_{|J|+\ell}$ on a neighborhood of the curve~$C_0$, and {\em smoothing parameters}, describing the ways to smooth out the nodes where $C_0$ meets the rest of the curve. Changing the moduli parameters always keeps $f$ in the stratum $P\ocX_{\l_1,\l_2,\dots,\l_r}$. The smoothing parameters were described in~\cite{ELSV2}. They are $u_i$, $1 \leq i \leq \ell$ and $a_{ij}$, $1 \leq j \leq \ell_i-1$ if we write the stable map in the following form:
\begin{equation} \label{Eq:HurwitzCoords}
f(z) = \sum_{i=1}^\ell \left[
\left(\frac{u_i}{z-z_i}\right)^{k_i} +
a_{i,1}\left(\frac{u_i}{z-z_i}\right) ^{k_i-1} +
\dots +  a_{i,k_i-1} \left(\frac{u_i}{z-z_i}\right)
\right].
\end{equation}

Introduce a coordinate~$z$ on the central component~$C_0$ of~$f$. Denote by $\a_1, \dots, \a_\ell$ the values of the $z$-coordinate at the intersection points with the branches. Denote by $\b_j, j \in J$ the values of the coordinate $z$ at the supplementary marked points. Then all functions $g(z)$ with zeros and poles of prescribed orders at the points $\a_i$ and $\b_j$ are proportional to
$$
g(z) = \frac{\prod_{j \in J}(z-\b_j)^{\l_j}}{\prod_{i=1}^\ell (z-\a_i)^{\s_i}}.
$$
We rewrite $g$ in the form~\eqref{Eq:HurwitzCoords} and denote by ${\bar u}_i$, ${\bar a}_{i,j}$ the values of the coefficients $u_i$ and $a_{ij}$ in this expansion.

Let $K = \mbox{LCM}(\s_1, \dots, \s_\ell)$ and $r_i = K/\s_i$ for $1 \leq i \leq \ell$. Then the stratum $P\ocX_{\l_1,\dots,\l_r}$ is parametrized in coordinates $u_i$ and $a_{ij}$ as follows:
$$
u_i = \zeta_i c^{r_i} \bar u_i, \quad a_{ij} = \zeta_i^j c^{j r_i} \bar a_{ij}.
$$
Here $(\zeta_1, \dots, \zeta_\ell)$ is a collection of roots of unity, $\zeta_i^{k_i} =1$. The action of $\Z/K \Z$ on $\Z/\s_1 \Z \times \dots \times \Z/ \s_\ell \Z$ has $\s_1 \cdots \s_\ell/K$ orbits corresponding to local irreducible components of the stratum $P\ocX_{\l_1,\dots,\l_r}$, and we choose one collection of roots of unity in each orbit. In each component of the stratum, the $\lambda_1$th derivative at~$x_1$ of the function corresponding to parameter~$c$ is equal to $c^K \cdot g^{(\lambda_1)}(x_1)$. Thus the vanishing order of the $\lambda_1$th derivative at $c=0$ equals $K$ for each component. Since there are $\s_1 \cdots \s_\ell /K$ components, we get the total vanishing order equal to $\s_1 \cdots \s_\ell$.

\bigskip

\noindent
{\bf e)}
{\em The vanishing order of $f^{(\l_1)}(x_1)$ along $P\ocX_{\l_1+1,\l_2\dots,\l_r}$ equals~$1$.} \\
This is similar to (d), but much easier. Once again, choose a generic $f \in P\ocX_{\l_1+1,\l_2\dots,\l_r}$ and a chart in  $\ocH_{r|\k}$, containing~$f$, that is a product of standard deformations of ramification loci. This time the important ramification locus is the supplementary marked point $x_i$. The corresponding deformation space is
$$
z^{\l_1+1} + a_1 z^{\l_1-1} + \dots + a_{\l_1}.
$$
The stratum $P\ocX_{\l_1,\l_2\dots,\l_r}$ in this coordinate is parametrized by one variable~$a$ as follows:
$$
(z-a)^{\l_1} (z+\l_1 a).
$$
The derivative of order $\l_1$ at $a$ is just $(\l_1+1)a$, thus its vanishing order along $a=0$ equals~1.

\bigskip

\noindent
{\bf f)} {\em The section $f^{(\l_1)}(x_1)$ has no other zeros.} \\
Let $f$ be any stable map lying in the closed stratum $P\ocX_{\l_1,\l_2\dots,\l_r}$ and satisfying $f^{(\l_1+1)}(x_1) = 0$. Denote by $X$ the connected component of $f^{-1}(0)$ containing~$x_1$. Let $J$ be the set of marked point contained in~$X$, $\ell$ the number of points at which $X$ meets other components of the curve and $\s'_1, \dots, s'_\ell$ the orders of zeros of $f$ at these points. If $X = \{x_1\}$ then we let $J = \{1 \}$, $\ell = 1$, $\s'_1 =$ order of the zero of $f$ at~$x_1$.

It follows from Lemma~\ref{Lem:closure} that all the connected components of $f^{-1}(0)$, including~$X$, are acceptable. In particular, we have $\sum_{j \in J} \l_j \leq \sum \s'_i$. Choose any list of integers $\s_i \leq \s'_i$ such that $\sum_{j \in J} \l_j = \sum \s'_i$. If $X = \{x_1\}$ we let $\s_1 = \l_1+1$. We claim that $f$ lies in the closure of the stratum $P\cX_{\l_1,\l_2,\dots,\l_r}^{J; \s_1, \dots, \s_\ell}$ (or of $P\cX_{\l_1+1,\l_2\dots,\l_r}$, if $X = \{ x_1 \}$). To see that apply to $f$ the following sequence of deformations. (i)~If $X$ is not smooth, deform it into a smooth curve. (ii)~Deform generically all ramification loci that do not contain supplementary marked points $x_i$ with $\l_i > 0$. (iii)~Deform all ramification loci over~$0$ except $X$ so that each supplementary marked point $x_i$ contained in these loci becomes a zero of order $\l_i$. This is possible by Lemma~\ref{Lem:closure}. (iv)~Deform the function $f$ at the $\ell$ attachment points to $X$ using the polynomial deformation families of degree~$\s'_i$ to reduce the orders of zeros at these points from $\s'_i$ to~$s_i$. Thus we have constructed a stable map $\tilde f \in P\ocX_{\l_1,\l_2\dots,\l_r}$ arbitrarily close to~$f$ that lies in the stratum $P\cX_{\l_1,\l_2,\dots,\l_r}^{J; \s_1, \dots, \s_\ell}$ or in $P\cX_{\l_1+1,\l_2\dots,\l_r}$. It follows that $f$ itself lies in the closure of one of these above strata. \QED

\subsection{Proof of the recursion}

The recursion formula of Theorem~\ref{Thm:Rec} is obtained by multiplying both sides of the cohomological relation~\eqref{Irel} by a monomial in the $\psi$-classes and taking the degrees of the corresponding cohomology classes.

\begin{lemma} \label{Lem:product}
The stratum $P\ocX_{\l_1,\l_2,\dots,\l_r}^{J; \s_1, \dots, \s_\ell}$ decomposes into several irreducible components each of which is a direct product of the moduli space $\ocM_{|J|+\ell}$ and $\ell$ multisingularity strata in Hurwitz spaces. Such components are in a one-to-one correspondence with the ways to split the set of poles and the set of supplementary marked points outside $J$ into $\ell$ unordered parts in such a way that each part contains at least one pole.
\end{lemma}

{\em Proof.} If we remove the central component from a curve in the stratum $P\ocX_{\l_1,\l_2,\dots,\l_r}^{J; \s_1, \dots, \s_\ell}$, the remaining part of the curve decomposes into $\ell$ connected components. Each of them contains at least one pole of the rational function~$f$, since otherwise $f$ would vanish identically on the connected component in question. The poles and the supplementary marked points lying on these connected components determine a splitting of the set of poles and supplementary marked points into $\ell$ parts.

Pick one splitting like that. Denote by $K_i$, $1 \leq i \leq \ell$, the list of orders of the poles on the $i$th connected component. Similarly, denote by $J_i$, $1 \leq i \leq \ell$, the set of supplementary marked points that lie on the $i$th connected component. Each element $s$ of $J_i$ comes with the corresponding label $\lambda_s$. Now, to each set $J_i$ we add one extra element corresponding to the point where the $i$th connected component meets the central component. This element carries the label~$\s_i$. Thus for each~$i$ we get the complete set of data that determines a multisingularity stratum in a Hurwitz space. Given $\ell$ functions in these $\ell$ Hurwitz spaces and one stable curve in $\ocM_{|J|+\ell}$ we can assemble a rational function $f$ from the intersection $P\ocX_{\l_1,\l_2,\dots,\l_r}^{J; \s_1, \dots, \s_\ell}$ by gluing them together in the natural way. Thus we obtain one irreducible component of the stratum $P\ocX_{\l_1,\l_2,\dots,\l_r}^{J; \s_1, \dots, \s_\ell}$.
\QED

\subsection{Proof of the recursion} \label{Ssec:ProofRec}

Recall that we must prove the following differential equation holds for the generating series $X$:
$$
\frac{\p X}{\p t_{s+1,m}}=\frac{\p X}{\p t_{s,m}}+s\,\frac{\p X}{\p t_{s,m+1}}-\sum_{\ell\ge1,a\ge0}\frac1{\ell!}\frac{\p\Psi_{a,\ell}}{\p t_{s,m}}
\sum_{\s_1+\dots+\s_\ell=a}\prod_{i=1}^\ell\Bigl(\s_i\,\frac{\p X}{\p t_{\s_i,0}}\Bigr).
$$

Let's consider the coefficient of the monomial
$$
\frac{\prod_{i=2}^r t_{s_i,m_i}}{|\Aut{\{ (s_i, m_i)\}}|}
$$
in this equality. The equality of these coefficients is obtained by multipying the cohomological identity~\eqref{Irel}
by $\psi_1^m \psi_2^{m_2} \dots \psi_r^{m_r}$ and taking the degree.

The coefficient in the left-hand side is the degree of
$$
[\ocX_{s+1, s_2, \dots, s_r}] \psi_1^m \psi_2^{m_2} \dots \psi_r^{m_r},
$$
which is the left-hand side of~\eqref{Irel} multiplied by $\psi_1^m \psi_2^{m_2} \dots \psi_r^{m_r}$.

The sum of the coefficients in the first two terms of the right-hand side is the degree of
$$
[\ocX_{s, s_2, \dots, s_r}] (\xi + s \psi_1) \psi_1^m \psi_2^{m_2} \dots \psi_r^{m_r},
$$
which is the first left-hand side term of~\eqref{Irel} multiplied by $\psi_1^m \psi_2^{m_2} \dots \psi_r^{m_r}$.

Finally, the terms of the last sum are indexed by a choice of $\ell$, a choice of $a$, a choice of $\s_1, \dots, \s_\ell$, and a way to decompose the monomial $t_{s_2, m_2} \cdots t_{s_k,m_k}$ into a product of the $\ell+1$ monomials $J, J_1, J_2, \dots, J_\ell$. The choices of $\ell$, $a$, $\s_1, \dots, \s_\ell$, and $J$ determine one term
$P\ocX_{\l_1,\l_2,\dots,\l_r}^{J; \s_1, \dots, \s_\ell}$
in the summation of the cohomological identity~\eqref{Irel}, where $a= \sum_{i \in J} \lambda_i$.  The choice of $J_1, \dots, J_\ell$
determines the choice of an irreducible component of the stratum $P\ocX_{\l_1,\l_2,\dots,\l_r}^{J; \s_1, \dots, \s_\ell}$,
as explained in Lemma~\ref{Lem:product}. Let us inspect the contribution to the coefficient of
$$
\frac{\prod_{i=2}^r t_{s_i,m_i}}{|\Aut{\{ (s_i, m_i)\}}|}
$$
arising from each given choice of $\ell$, $a$, $\s_1, \dots, \s_\ell$, $J, J_1, \dots, J_\ell$. This contribution is a product of several factors.

The first factor is
$$
\frac{|\Aut\{ (s_i, m_i)\}|}{|\Aut(J)| \prod_{i=1}^{\ell} |\Aut(J_i)|}.
$$
This factor accounts for the fact that the supplementary marked point are numbered, while the corresponding variables $t_{s_i,m_i}$ in the monomial are not. The factor transforms the number of ways to split the monomial into a product into the number of ways to split the supplementary marked points into $\ell+1$ groups.

The second factor is $1/\ell!$. It accounts for the fact that the irreducible components of the stratum $P\ocX_{\l_1,\l_2,\dots,\l_r}^{J; \s_1, \dots, \s_\ell}$
are labelled by unordered set partitions into $\ell$ parts, while the terms of the product
$$
\prod_{i=1}^\ell \left( \s_i \frac{\partial X}{\partial t_{\s_i,0}} \right)
$$
are ordered.

The third factor is $\prod_{i=1}^\ell \s_i$. This factor is present in the cohomological identity~\eqref{Irel}, where it represents the vanishing order of the section.

The fourth factor is the coefficient of the series $\Psi_{a,\ell}$. As explained in Remark~\ref{Rem:Psi}, this coefficient is the intersection number
$$
\int_{\ocM_{|J|+\ell}} \prod_{i \in J} \psi_i^{m_i}.
$$
This factor appears as we integrate over the moduli of the central component the product of relevant $\psi$-classes.

As stated in Lemma~\ref{Lem:product}, each irreducible component of the stratum $P\ocX_{\l_1,\l_2,\dots,\l_r}^{J; \s_1, \dots, \s_\ell}$ is the product of the moduli space $\ocM_{|J|+\ell}$ and $\ell$ multisingularity strata in Hurwitz spaces. Moreover, each multisingularity stratum comes with a monomial in $\psi$-classes. So as not to multiply unnecessary notation, let's just denote the multisingularity strata by $\ocX^{(i)}$ and the corresponding monomials in $\psi$-classes by $\widetilde{\psi}^{(i)}$ for $1 \leq i \leq \ell$. Then the fifth and last factor is the product of degrees
$$
\prod_{i=1}^\ell \deg \,
[\ocX^{(i)} \widetilde{\psi}^{(i)}].
$$

Now it remains to explain why the product of factor~4 (the integral over the moduli of the central component) and factor~5 (the product of degrees of $\ell$ multisingularity strata multiplied by monomials $\widetilde{\psi}^{(i)}$) is equal to the degree of the corresponding irreducible component of the stratum $P\ocX_{\l_1,\l_2,\dots,\l_r}^{J; \s_1, \dots, \s_\ell}$ multiplied by $\prod \psi_i^{m_i}$.

This follows from a general statement. The Segre class of a product of cones $C_i \to B_i$ is equal to the product of their Segre classes. Given a collection of cohomology classes $\a_i \in H^*(B_i)$, the degree of the product $\prod_i \a_i$ is equal to the evaluation on $\prod B_i$ of $\prod \a_i \cdot s(\prod C_i)$. The latter, in turn, is equal to the product of evaluations of $\a_i \cdot s(C_i)$ on $B_i$. In our case, the cones $C_i$ are the multisingularity strata $\ocX^{(i)}$ and one trivial one over $\ocM_{|J|+\ell}$. The classes $\a_i$ are the monomials in $\psi$-classes.

We conclude that the recursion equality is obtained by multiplying the cohomological relation~\ref{Irel} by the monomial $\psi_1^{m_1} \cdots \psi_r^{m_r}$ and taking the degrees of all classes. This completes the proof of Theorem~\ref{Thm:Rec}. \QED

\section{A reduced form of the recursion}

Let us set $p_i=t_{i,0}$ $(i\ge1)$ and $t_j=t_{0,j}$ $(j\ge 0)$.  The variables $p_1,p_2,\dots;t_0,t_1,\dots$ form the \emph{reduced set of the $t$-variables}. Consider the restriction homomorphism
\begin{gather*}
\Q[t_{0,0},t_{1,0},t_{0,1},\dots]\longrightarrow\Q[p_1,p_2,\dots;t_0,t_1,\dots],\\
t_{i,0}\mapsto p_i~(i\ge1),\quad t_{0,j}\mapsto t_j,\quad t_{i,j}\mapsto 0~(i,j\ge1).
\end{gather*}
Denote by $\oX$ and $\oX_{s,m}$ the image of the descendant Hurwitz potential~$X$ and its partial derivative $\p_{t_{s,m}}X$, respectively, under the reduction. Then one can observe that the collection of functions~$\oX$ and~$\oX_{s,m}$, $s,m\ge0$, is closed with respect to the recursion of Theorem~\ref{Thm:Rec}.

\begin{theorem}\label{thred}
The functions~$\oX_{s,m}$ considered as formal power series in the variables $p_1,p_2,\dots;t_0,t_1,\dots$ obey the following equations
\begin{align*}
\p_{p_s}\oX&=\oX_{s,0}\quad(s\ge1),\\
\oX_{s+1,m}&=\oX_{s,m}+s\,\oX_{s,m+1}
- \sum_{\ell\ge1,~a\ge s}\frac1{\ell!}\frac{\p\oPsi_{a-s,\ell}}{\p t_{m}}
\sum_{\s_1+\dots+\s_\ell=a}\prod_{i=1}^\ell\bigl(\s_i \oX_{\s_i,0}\bigr)\quad(m,s\ge0),\\
\oX_{0,m}&=\p_{t_m}\oX\quad(m\ge0),
\end{align*}
where
$$\oPsi_{a,\ell}=\sum_{k,j\ge0}
\left(\frac1{j!} \sum_{\substack{\l_1+\dots+\l_j=a\\\l_i\ge1}}
  p_{\l_1}\dots p_{\l_j}\right)
\left(\frac1{k!}
 \sum_{\substack{\nu_1+\dots+\nu_k=\ell+k+j-3\\\nu_i\ge0}}
 \binom{|\nu|}{\nu_1,\dots,\nu_k}t_{\nu_1}\dots t_{\nu_k}\right)
$$
\end{theorem}

Using the reduced set of variables means that, for the Hurwitz spaces $\ocH_{r|\k}$ under consideration,
 we associate both degeneracy conditions and the $\psi$-classes to only one of the supplementary marked points,
 say, to the first one, while to the other supplementary marked points we associate either
 degeneracy conditions or the $\psi$-classes, but not both.

Applying relations of this theorem repeatedly, one expresses each partial derivative $\p_{p_s}\oX=\oX_{s,0}$ in the direction of the $p$-variables as a power series in the variables~$t_i$, $p_j$, and the partial derivatives $\p_{t_m}\oX=\oX_{0,m}$ in the direction of the $t$-variables. Therefore, these relations determine the series $\oX$ uniquely from the initial conditions
$$\oX\Bigm|_{p_j=0}=\sum_{r\ge0}\frac1{r!}
  \sum_{\nu_1,\dots,\nu_r}
 \binom{|\nu|}{\nu_1,\dots,\nu_r}z_{|\nu|,r}\;t_{\nu_1}\dots t_{\nu_r},$$
where $z_{d,r}=z_{d,r}(q)$ is the same as in Proposition~\ref{Prop:Init}.

Theorem~\ref{thred} provides a reduced and simplified form of the recursion sufficient for an independent computation of the Hurwitz numbers of Proposition~\ref{prophlambda}. Some further reduction is discussed in the next section.

\section{String and dilaton equations}

A study of the forgetful map $\pi:\ocH_{r+1|\k}\to\ocH_{r|\k}$ given by forgetting the last supplementary marked point leads to partial differential equations called the \emph{string} and \emph{dilaton equations}. They express the partial derivatives~$\p_{t_{0,0}}X$ and~$\p_{t_{0,1}}X$, respectively, in terms of partial derivatives in other variables.

\begin{theorem} The descendant genus zero Hurwitz potential $X(t,q)$ obeys the following string and dilaton equations
\begin{align*}
\frac{\p X}{\p t_{0,0}}&=\sum_{\l,d\ge0}t_{\l,d+1}\,\frac{\p X}{\p t_{\l,d}}
 +\sum_{k\ge1}k\,q_k\frac{\p X}{\p q_k},\\
\frac{\p X}{\p t_{0,1}}&=\sum_{\l,d\ge0}t_{\l,d}\,\frac{\p X}{\p t_{\l,d}}
 +\sum_{k\ge1}q_k\frac{\p X}{\p q_k}-2\,X.
\end{align*}
\end{theorem}

\paragraph{Proof.} Denote the marked poles of $f$ by $p_1, \dots, p_n$, and the supplementary marked points by $x_1, \dots, x_r$ or $x_1, \dots, x_{r+1}$.

Consider the locus $\delta_i \subset \ocH_{r+1|\k}$ composed by stable maps $f: C \to \CP^1$ of the following form. The source curve $C$ has $n$ disjoint irreducible components $C_1$, \dots, $C_n$
containing the marked points $p_1, \dots, p_n$ respectively; the restriction of $f$ to $C_1$, \dots, $C_n$ has the form $z \mapsto z^{k_1}$, \dots, $z^{k_n}$; these $n$ components are connected by a union $X$ of one or several irreducible components on which the map~$f$ vanishes identically;
 the subcurve $X$ contains the supplementary marked points $x_1, \dots, x_r$; the marked point $x_{r+1}$ belongs to the component $C_i$.
\begin{center}
\includegraphics[width=10em]{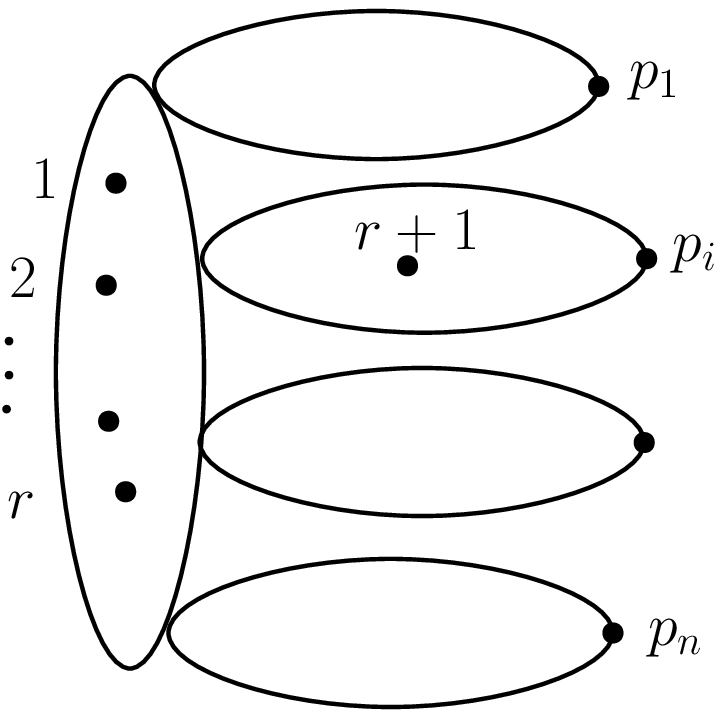}
\end{center}

The reason we have introduced the loci $\delta_i$ is that the forgetful map $\pi$ maps them to the zero section of the cone $\ocH_{r|\k}$. Thus the forgetful map $P\ocH_{r+1|\k}\to P\ocH_{r|\k}$ is not defined on $P\delta_i$. To construct a well-defined forgetful map we introduce the blow-up $\Bl$ of $P\ocH_{r+1|\k}$ along the loci $P\delta_i$, and denote by $\Delta_i$ the exceptional divisor corresponding to $P\delta_i$. Thus we obtain the following commutative diagram.
\begin{center}
\includegraphics[width=17em]{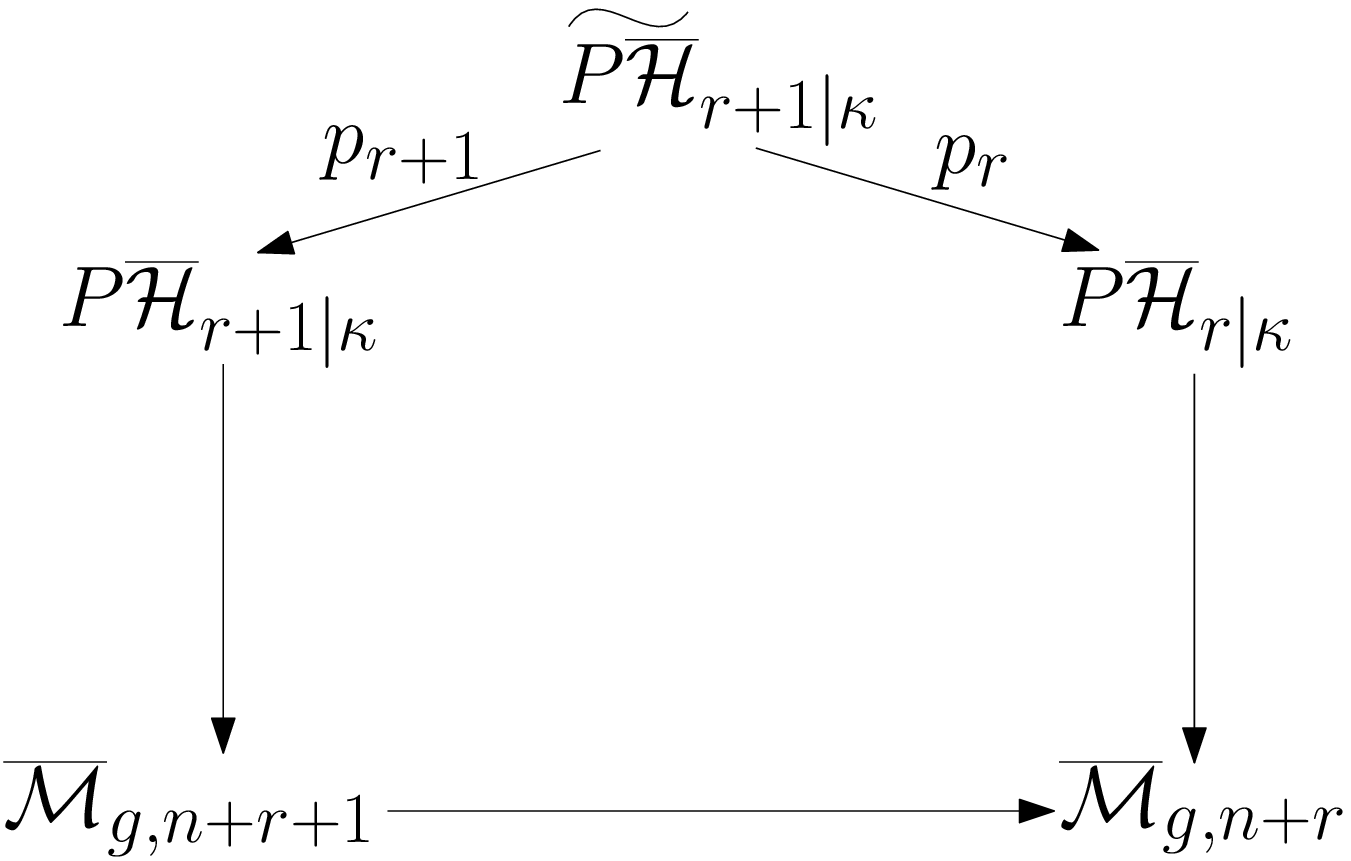}
\end{center}

Note that the locus $\delta_i$ is isomorphic to the quotient orbifold
$$
\ocM_{0, r+n} / \left( \Z / k_1 \Z \times \dots \times \widehat{\Z /k_i \Z} \times \dots \times \Z / k_n \Z \right)
$$
with the trivial group action, where the hat means a skipped factor. The restriction of $p_r$ to the corresponding exceptional divisor $\Delta_i$ is a degree~$k_i$ map onto $P\ocH_{r|\k}$.

We will work in the blown-up space $\Bl$. Denote by $\psi_1, \dots, \psi_{r+1}$ the pull-backs of the $\psi$-classes from $\ocM_{g,n+r+1}$ at the supplementary marked points to $\Bl$. Similarly, denote by $\opsi_1, \dots, \opsi_r$ the pull-backs of the $\psi$-classes from $\ocM_{g,n+r}$ at the supplementary marked points to $\Bl$. Further, denote by $\xi$ the pull-back of the $\xi$-class from $P\ocH_{r+1|\k}$ to $\Bl$ and by $\oxi$ the pull-back of the $\xi$-class from $P\ocH_{r|\k}$ to $\Bl$. Finally, denote by $D_j$, $1 \leq j \leq r$, the pull-back from $\ocM_{g,n+r+1}$ of the divisor of curves of the form
\begin{center}
$D_j = \Biggl\{ \includegraphics[width=15em,trim=0 3em 0 0]{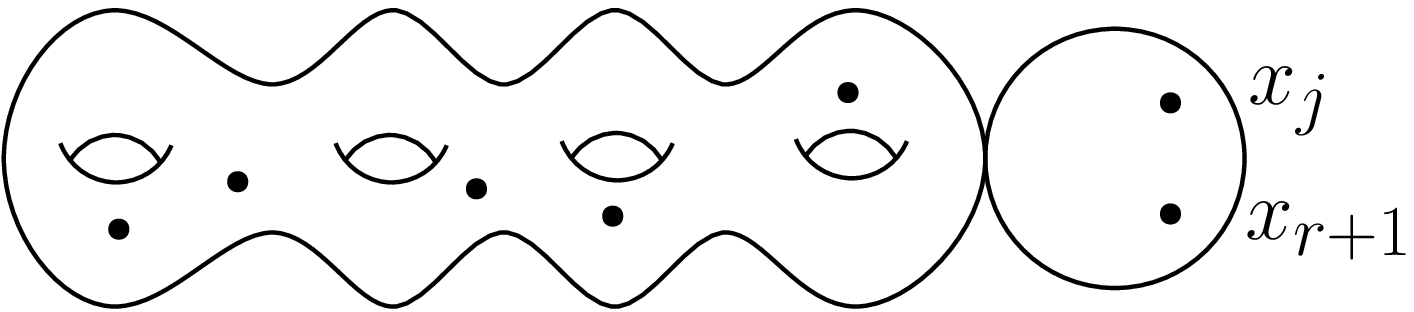} \Biggr\}$.
\end{center}
Note that the restriction of $p_r$ to $D_j$ is an isomoprhism onto $P\ocH_{r|\k}$.

\begin{lemma} \label{Lem:relations}
The following relations hold in $H^*(\Bl)$:
\begin{itemize}
\item[a)]
$D_j D_{j'} = 0$ for $j \not= j'$;
\item[b)]
$\Delta_i D_j = 0$ for $1 \leq i \leq n$, $1 \leq j \leq r$;
\item[c)]
$\Delta_i \Delta_{i'} = 0$ for $i \not= i'$;
\item[d)] $\xi \Delta_i =0$ for $1 \leq i \leq n$;
\item[e)] $\psi_j \D_j=0$ for $1 \leq j \leq r$;
\item[f)] $\psi_{r+1} \D_j=0$ for $1 \leq j \leq r$;
\item[g)] $\psi_{r+1} \Delta_i=0$ for $1 \leq i \leq n$.
\item[h)]
$\xi = \oxi + \sum_{i=1}^n \Delta_i$;
\item[i)]
$\psi_j = \opsi_j + D_j$ for $1 \leq j \leq r$;
\end{itemize}
\end{lemma}

\paragraph{Proof.} Relations (e), (f), and~(i) are pull-backs of the well-known analogous relations from $\ocM_{g,n+r+1}$. Relations (a), (b), and~(c) follow from the fact that the geometric intersections of the corresponding divisors are empty. Relations (d) and~(g) express the vanishing of the intersection between the exceptional divisor of a blow-up and a pull-back class from the base. Relation~(h) is slightly harder. The line bundles $p_{r+1}^*(\cO(1))$ and $p_r^*(\cO(1))$ are naturally identified except over the exceptional divisors $\Delta_i$. So we have $\xi = \oxi + \sum c_i \Delta_i$ for some coefficients~$c_i$. Now, the normal line bundle to the exceptional divisor $\Delta_i$ is identified with $\cO(-1)$. Thus we have the equation $\xi \Delta_i =0 = \oxi \Delta_i - c_i \oxi \Delta_i$. It follows that $c_i = 1$ for every~$i$.
\QED

It follows from the relations that
$$
\psi_j^d = \opsi_j^d + D_j \opsi^{d-1},
$$
$$
\xi_j^d = \oxi_j^d + \sum_{i=1}^n \Delta_i \oxi^{d-1}.
$$
This can be easily proved by induction on~$d$.

Now, consider a stratum $P\ocX_{\l_1,\l_2,\dots,\l_r} \subset P\ocH_{r|\k}$, and the stratum $P\ocX_{\l_1,\l_2,\dots,\l_r,0} \subset P\ocH_{r+1|\k}$ obtained by adding one supplementary marked point with $\l_{r+1} =0$. The preimage $p_r^{-1}(P\ocX_{\l_1,\l_2,\dots,\l_r})$ is the proper transform of $P\ocX_{\l_1,\l_2,\dots,\l_r,0}$ under the blow-up $p_{r+1}$. For shortness, we will denote the proper transform by $\tX$ and the strata by $P\ocX_r$ and $P\ocX_{r+1}$, respectively.

Now we can prove the dilaton relation. The integral
$$
\int_{P\ocX_{r+1}} \xi^d \psi_1^{d_1} \dots \psi_r^{d_r} \psi_{r+1},
$$
where $d$ is determined by degree reasons, is the coefficient of
$$
\frac{\prod_{j=1}^r t_{\l_j, d_j} \prod_{i=1}^n q_{k_i}}
{|\Aut \{ (\l_j, d_j)\} | \cdot |\Aut \{ k_i\} |  }
$$
in the power series~$\p X / \p t_{0,1}$. Applying the projection formula to $p_{r+1}$ we rewrite the integral as
$$
\int_{\tX} \xi^d \psi_1^{d_1} \dots \psi_r^{d_r} \psi_{r+1}.
$$
Now we use the relations from Lemma~\ref{Lem:relations}:
$$
\int_{P\ocX_{r+1}} \xi^d \psi_1^{d_1} \dots \psi_r^{d_r} \psi_{r+1} =
$$
$$
\int_{\tX} \left(\oxi^d + \sum k_i \Delta_i \oxi^{d-1}\right)
(\opsi_1^{d_1} + D_1 \opsi_1^{d_1-1})   \dots
(\opsi_r^{d_r} + D_r \opsi_r^{d_r-1})
\psi_{r+1}=
$$
$$
\int_{\tX} \oxi^d \opsi_1^{d_1}
\opsi_r^{d_r}
\psi_{r+1}.
$$
Applying the projection formula to the map $p_r$, we get
$$
(n+r-2) \int_{P\ocX_r} \oxi^d \opsi_1^{d_1}
\opsi_r^{d_r}.
$$
This is the right-hand side of the dilaton relation, where the multiplication of the coefficient of $X$ by $n+r-2$ is realized by the differential operator
$$
\sum_{\l,d\ge0}t_{\l,d}\,\frac{\p}{\p t_{\l,d}}
 +\sum_{k\ge1}q_k\frac{\p }{\p q_k}-2.
$$
Thus we have proved the dilaton relation.

We now prove the string relation in a similar way. The integral
$$
\int_{P\ocX_{r+1}} \xi^d \psi_1^{d_1} \dots \psi_r^{d_r},
$$
where $d$ is determined by the degree argument, is the coefficient of
$$
\frac{\prod_{j=1}^r t_{\l_j, d_j} \prod_{i=1}^n q_{k_i}}
{|\Aut \{ (\l_j, d_j)\} | \cdot |\Aut \{ k_i\} |  }
$$
in the power series~$\p X / \p t_{0,0}$. Applying the projection formula to $p_{r+1}$ we rewrite the integral as
$$
\int_{\tX} \xi^d \psi_1^{d_1} \dots \psi_r^{d_r}.
$$
Now we use the relations from Lemma~\ref{Lem:relations}:
$$
\int_{P\ocX_{r+1}} \xi^d \psi_1^{d_1} \dots \psi_r^{d_r} =
$$
$$
\int_{\tX} \left(\oxi^d + \sum k_i \Delta_i \oxi^{d-1}\right)
(\opsi_1^{d_1} + D_1 \opsi_1^{d_1-1})   \dots
(\opsi_r^{d_r} + D_r \opsi_r^{d_r-1})
=
$$
$$
\sum_{i=1}^n \int_{\Delta_i} \oxi^{d-1} \opsi_1^{d_1} \dots \opsi_r^{d_r}
+ \sum_{j=1}^r \int_{D_j} \oxi^d \opsi_1^{d_1}\dots
 \opsi_j^{d_j-1}  \dots  \opsi_r^{d_r}.
$$
Applying the projection formula to the map $p_r$, and recalling that $p_r|_{D_j}$ has degree~1, while $p_r|_{\Delta_i}$ has degree~$k_i$, we get
$$
 \sum_{i=1}^n k_i  \int_{P\ocX_r}\oxi^{d-1} \opsi_1^{d_1} \dots \opsi_r^{d_r} + \sum_{j=1}^r  \int_{P\ocX_r} \oxi^d \opsi_1^{d_1}\dots
 \opsi_j^{d_j-1}  \dots  \opsi_r^{d_r}.
$$
This is the right-hand side of the string relation. \QED

The string and dilaton equations provide a further reduction of the recursion of Theorem~\ref{thred}. Namely, they allow one to eliminate the variables~$t_0=t_{0,0}$ and~$t_1=t_{0,1}$ from consideration by setting them to be equal to zero.

Note that the $q$-derivatives of the right-hand sides of these equations preserve the ring generated by the functions~$z_{d,r}(q)$ of Proposition~\ref{Prop:Init}. Namely, their action on the functions~$z_{d,r}$ is given by the following equations
\begin{equation}
\begin{aligned}
\Bigl(\sum_{k\ge1}k\,q_k\frac{\p }{\p q_k}\Bigr)z_{d,r}&=z_{d,r+1}-z_{d-1,r},\\
\Bigl(\sum_{k\ge1}q_k\frac{\p }{\p q_k}\Bigr)z_{d,r}&=(d+1)\,z_{d+1,r+1}+(2-r)\,z_{d,r}.
\end{aligned}\label{eqzred}\end{equation}
These relations allow one to simplify manipulations with the potentials $X$ and $\oX$ by keeping the Taylor coefficients of their $t$-expansions as polynomials in the functions~$z_{d,r}(q)$.

The proof of~\eqref{eqzred} is obvious: these relations are equivalent to the following identities between the binomial coefficients:
\begin{align*}\binom{n+r-3}{d}&=\binom{n+r+1-3}{d}-\binom{n+r-3}{d-1},\\
n\,\binom{n+r-3}{d}&=(d+1)\,\binom{n+r+1-3}{d+1}+(2-r)\,\binom{n+r-3}{d}.
\end{align*}

We conclude this section with a simple remark. Since the coefficients of the series $\Psi_{a,\ell}$ of Remark~\ref{Rem:Psi} are intersection numbers of $\psi$-classes over genus~0 moduli spaces, these series also satisfy a string and a dilaton equations.

\begin{proposition} The series $\Psi_{a,\ell}$ of Theorem~\ref{Thm:Rec} obey the following string and dilaton equations
\begin{align*}
\frac{\p\Psi_{a,\ell}}{\p{t_{0,0}}}&=\Psi_{a,\ell+1}
 +\sum_{i,j\ge0}t_{i,j+1}\frac{\p\Psi_{a,\ell}}{\p{t_{i,j}}},\\
\frac{\p\Psi_{a,\ell}}{\p{t_{0,1}}}&
 =\sum_{i,j\ge0}t_{i,j}\frac{\p\Psi_{a,\ell}}{\p{t_{i,j}}}+(\ell-2)\,\Psi_{a,\ell}.
\end{align*}
\end{proposition}

\section{Computing residual polynomials}

The goal of this section is to relate the statements of Theorems~\ref{Thm:Rec} and~\ref{cohth} to the general multisingularity principle.

Let $F:X\to Y$ be a generic finite holomorphic mapping of two compact complex manifolds
of the same dimension. The mapping is {\it degenerate\/} at a point $x\in X$
if $dF(x)$ is a degenerate linear mapping of the tangent spaces at~$x$ and~$F(x)$,
respectively. For a given singularity type~$\sigma$, denote by $\Sigma_\sigma(X)$
the closure of the subset of points in~$X$ where~$F$ attains a singularity
of type~$\sigma$. The {\it Thom principle\/} states that the cohomology class
$[\Sigma_\sigma(X)]\in H^*(X)$ Poincar\'e dual to the subvariety
$\Sigma_\sigma(X)$ admits a universal expression
as a polynomial in relative Chern classes $c_0(F),c_1(F),\dots$
of the mapping~$F$ defined by the series expansion
$$
c(F)=\frac{c(F^*(TY))}{c(TX)}=c_0(F)+c_1(F)+\dots, \qquad c_i(F)\in H^{2i}(X).
$$
The polynomial is quasihomogenous, of degree equal to the codimension of~$\Sigma_\sigma(X)$.

This principle has been extended to the case of multisingularities
by M.~Kazarian~\cite{K1}. Namely, let $\sigma_1,\dots,\sigma_m$ be a set of
singularity types. Denote by $\Sigma_{\sigma_1,\dots,\sigma_m}(Y)$
the closure of the locus of points in~$Y$, whose preimages
contain singularities of types $\sigma_1,
\dots,\sigma_m$. Let $[\Sigma_{\sigma_1,\dots,\sigma_m}(Y)]\in H^*(Y)$
be the corresponding Poincar\'e dual cohomology class.
Then there are universal polynomials
$R_{\sigma_1,\dots,\sigma_m}$ in the classes~$c_i(F)$
such that the generating function
$\sum[\Sigma_{\sigma_1,\dots,\sigma_m}(Y)]t_{\sigma_1}\dots t_{\sigma_m}$,
where the summation is carried over all tuples of singularity types,
is the exponent $\exp~F_*\cR$ of  the pushforward, under~$F$, of the generating function
$\cR=\sum R_{\sigma_1,\dots,\sigma_m}t_{\sigma_1}\dots t_{\sigma_m}$.

We shall apply the latter principle to the case of families of meromorphic functions;
its applicability in this case has been proved in~\cite{KL1a}. Let $\cH=\ocH_\k$ be one of the Hurwitz spaces studied in this paper. The points of this space represent certain meromorphic functions. This family of functions fits into a diagram of mappings
$$
\xymatrix {\widetilde{X}\ar[rr]^{\widetilde F}\ar[rd]_{\widetilde P}&&
\widetilde{Y}\ar[ld]^{\widetilde {Q}}\\&\cH}
$$
Here $\widetilde{X}$ is the universal curve over $\cH$ that can be identified with $\ocH_{1|\k}$ in the notation of the present paper,
and $\widetilde{Y}=\cH\times \C P^1$. The multiplicative group $\C^*$ of nonzero complex numbers acts on the space~$\cH$
by multiplying functions by constants. This action extends naturally to the
spaces~$\widetilde X$ and $\widetilde{Y}$. It preserves the singularities and their types.
Denote the $\C^*$-quotient spaces by $P\cH$, $X$, and $Y$, respectively, and the
corresponding quotient mappings by $P:X\to P\cH$, $Q:Y\to\cH$, $F:X\to Y$,
so that we have a commutative diagram of the form
$$
\xymatrix {X\ar[rr]^{F}\ar[rd]_{P}&&
{Y}\ar[ld]^{{Q}}\\&P\cH}
$$
We say that~$\widetilde{F}$ {\it has a singularity of type~$A_{k-1}$ at a point~$x\in \widetilde{X}$}
if~$x$ is a nonsingular point of the fiber of~$\widetilde{P}$ passing through~$x$,
and the restriction of~$\widetilde{F}$ to this fiber admits, in an appropriate local
coordinate~$z$ in the fiber of~$\widetilde{P}$ at~$x$ and an appropriate local
coordinate in the fiber of~$\widetilde{Q}$ at~$\widetilde{F}(x)$ the form $z\mapsto z^k$.
Since multiplication of a function by a nonzero constant does not change the type of a singularity at a point,
the same definition remains valid for a point $x\in X$.
The shift by~$1$ in the notation of the singularity type is due to the tradition coming from singularity theory. $A_{k-1}$ singularities
are the only possible isolated singularity types at smooth points of the fibers of~$P$.
Since the singularity types~$A_{k-1}$, $k=1,2,3,\dots$ are indexed by integer numbers, the multisingularity classes are labeled by tuples of integers: the locus of functions with ramifications of orders $\nu_1,\dots,\nu_r$ over one point is denoted by $\Sigma_{\nu_1,\dots,\nu_r}(Y)$. Respectively, the generating function~$\cR$ for the corresponding residual polynomials
can be considered as a function in formal variables $t_1,t_2,\dots$,
where the variable~$t_k$ is in charge of the singularity type~$A_{k-1}$; in particular,
$t_1$ corresponds to nonsingular points of the mapping~$F$.

If the fibers of~$P$ are smooth, the universal polynomials are computed in~\cite{KL1a}. We reproduce these computations below. Unfortunately, these computations do not provide an answer in the general case when the fibers of~$P$ are allowed to be singular. Moreover, the whole multisingularity formula is not applicable in its straightforward form since the genericity conditions required for its applicability break down in the presence of nonisolated singularities. Thus, the formulas obtained in~\cite{KL1a} should be corrected by terms supported on the loci of nonisolated singularities.

Observe that the class of the multisingularity locus $\Sigma_{\nu_1,\dots,\nu_r}(Y)$ can also be computed by pushing forward the class of the subvariety $X(\nu_1,\dots,\nu_r)\subset P\ocH_{r|\k}$ discussed in the previous sections of the present paper under the natural forgetful map $P\ocH_{r|\k}\to P\ocH_{1|\k}\to P\ocH_{\k}$. In particular, our inductive computation of the classes $[X(\nu_1,\dots,\nu_r)]$ undertaken in this paper implies implicitly the computation of all necessary correction terms for the classes $R_{\nu_1,\dots,\nu_r}$ participating in the general multisingularity formula. However, the explicit computation of the pushforward homomorphism is still to be done.

We show below that the part of the series $\cR$ corresponding to the contribution of the smooth fibers of~$P$ satisfies the equations
of the Kadomtsev--Petviashvili (KP) integrable hierarchy of partial differential equations.
One can hope that the corresponding enriched generating functions involving the contribution of
nonisolated singularities also are solutions to suitable integrable hierarchies.

It is proved in~\cite{KL1a} that in the case where the fibers of~$P$ are smooth,
all the relative Chern classes of~$F$ can be expressed as universal polynomials in just two
classes, $\xi,\psi\in H^2(X)$. Namely, the class $\psi=\psi_1$ is the $\psi$ class associated with the unique supplementary marked point on $X=P\ocH_{1|\k}$. The class $-\psi$ can be defined also as the relative first Chern class of the fibration~$P:X\to P\cH$.
The class $\xi$ is the first Chern class~$c_1(\cO(1))$,
where the line bundle $\cO(1)$ is inherited from the $\C^*$-action on~$\widetilde{X}$.
Hence, we define
$$
\cR(\psi,\xi;t_1,t_2,\dots)=\sum R_{i_1,i_2,\dots}(\psi,\xi)t_{i_1}t_{i_2}\dots.
$$
The first few terms of the function~$\cR$ are
$$
\cR=t_1+\left(-\frac12 t_1^2+(\xi+\psi)t_2\right)
+\left(\frac13t_1^3-2(\xi+\psi)t_1t_2+(\xi+\psi)(\xi+2\psi)t_3\right)+\dots.
$$

\begin{theorem} For a generic family of functions on smooth curves,
the generating function~$\cR$ of residual polynomials is a solution to the scaled KP hierarchy of partial
differential equations. In particular, it solves the first scaled KP equation of the form
$$
\frac{\partial^2\cR}{\partial t_2^2}=2\psi\xi\left(\frac{\partial^2\cR}{\partial t_1^2}\right)^2
+\frac43\frac{\partial^2\cR}{\partial t_1\partial t_3}
-\frac13\psi^2\frac{\partial^4\cR}{\partial t_1^4}.
$$
\end{theorem}

The {\it scaled KP equations\/} are obtained from the ordinary ones by applying
to them the following scaling: a partial
derivation of order~$k$ is replaced by $\xi(-\psi)^{k-1}$ times the same derivation.
Since all the partial derivatives in the KP equations are of order at least~$2$,
the scaled equations are divisible by~$\xi\psi$, and we simplify them by dividing
by this monomial.

{\bf Proof.} All the solutions of the KP equations are known to be logarithms of
 tau functions for the KP-hierarchy. The tau function
 $\exp(-\frac\xi\psi\cR A)$ is, in fact, the following one:
 $$
 \exp(-\frac\xi\psi\cR A)=1+\xi\tilde s_1+\xi(\xi+\psi)\tilde s_2
 +\xi(\xi+\psi)(\xi+2\psi)\tilde s_3\dots,
 $$
 where $\tilde s_k$ is the $k$th {\it scaled one-part Schur polynomial\/} (written in the variables~$t_i$), which are the homogeneous parts of the decomposition
 $$
 e^{-\frac1\psi(t_1+t_2+t_3+\dots)}= \tilde s_0+\tilde s_1+\tilde s_2+\dots
 $$
 As is well known, an arbitrary linear combination of one-part Schur polynomials
 is a tau function.

 In order to check that these formulas indeed produce the correct universal
 expressions for the residual polynomials, it suffices to verify them on a sufficiently
 rich bunch of examples of Hurwitz spaces. For such a bunch, consisting of the
 spaces of polynomials (versal unfoldings of the singularities~$A_k$, $k=1,2,\dots$)
 the verification has been done in~\cite{KL1a}. \QED


\begin{thebibliography}{99}

\bibitem{ACG} E.~Albarello, M.~Cornalba, P.~Griffiths,
{\it Geometry of algebraic curves Vol.~II}, Springer (2011)

\bibitem{B} B.~Bychkov,
{\it Degrees of cohomology classes of multsingularities in Hurwitz spaces of rational functions},
Funct. Anal. Appl., {\bf 53}, no.~1, 11--22 (2019)

%
%
\bibitem{ELSV2} T.~Ekedahl, S.~K.~Lando, M.~Shapiro, A.~Vainshtein,
{\it Hurwitz numbers and intersections on moduli spaces of
curves}, Invent. math., {\bf 146}, 297--327~(2001)

\bibitem{Fulton} W.~Fulton. {\em Intersection theory}, vol.~2 of Ergebnisse der Mathematik und ihrer Grenzgebiete. 3. Folge. A Series of Modern Surveys in Mathematics [Results in Mathematics and
Related Areas. 3rd Series. A Series of Modern Surveys in Mathematics]. Springer-Verlag,
Berlin, second edition, 1998.

\bibitem{GJ} I.~P.~Goulden, D.~M.~Jackson,
{\it Transitive factorisation into transpositions and holomorphic
mappings on the sphere}, Proc. Amer. Math. Soc., {\bf 125}, no.~1,
51--60~(1997)

%

\bibitem{K1}  M.~Kazarian, {\it Multisingularities, cobordisms, and enumerative geometry},
Uspekhi Mat. Nauk, no 4, \textbf{58} (2003), 29--88; translation
in Russian Math. Surveys, (4) \textbf{58} (2003), 665-724.


\bibitem{KL1} M.~Kazarian, S.~Lando,
{\it On intersection theory on Hurwitz spaces}, Izv. Ross. Akad.
Nauk Ser. Mat., {\bf 68}, no. 5, 91--122~(2004); translation in
Izv. Math. {\bf 68}, no. 5, 935--964~(2004)

\bibitem{KL1a} M.~Kazarian, S.~Lando, {\it Thom polynomials for mappings of curves with
isolated singularities}, in Tr. Mat. Inst. Steklova (2007) Anal. i Osob. Ch. 1;
translated in: Proc. Steklov Inst. Math. (2007), no.~1, 93--106

%

\bibitem{KL4} M.~Kazarian, S.~Lando, {\it Combinatorial solutions to integrable
hirarchies},  Uspekhi Mat. Nauk 70 (2015), no. 3 (423), 77--106.
English translation: 2015 Russ. Math. Surv. vol.~70, 453--482

\bibitem{KLZ} M.~Kazarian, S.~Lando, D.~Zvonkine,
{\it Universal cohomological expressions for singularities in
families of genus~0 stable maps}, International Mathematics Research Notices 22 (2018), 6817--6843, \texttt{doi.org/10.1093/imrn/rnx070}.


\bibitem{LZ} S.~K.~Lando, A.~K.~Zvonkin,
{\it Graphs on surfaces and their applications}, Springer (2004)

\bibitem{LZ1} S.~K.~Lando, D.~Zvonkine,
{\it On multiplicities of the Lyashko-Looijenga mapping on strata
of the discriminant}, Funktsional. Anal. i Prilozhen., {\bf 33},
no. 3, 21--34, (1999); translation in Funct. Anal. Appl., {\bf
33}, no. 3, 178--188  (1999)

\bibitem{LZ2} S.~K.~Lando, D.~Zvonkine,
{\it Counting Ramified Coverings and Intersection Theory on Spaces
of Rational Functions I (Cohomology of Hurwitz Spaces)}, Moscow
Math.~J., {\bf 7} (1), 85--107~(2007)

\bibitem{Ok1} A.~Okounkov, {\it Toda equations for Hurwitz numbers},
Math. Res. Lett. 7 (2000), no.\ 4, 447--453.

\bibitem{Sauvaget} A.~Sauvaget, {\it Cohomology classes of strata of differentials.}, \texttt{arXiv:1701.07867}.

\bibitem{W1}
 E.~Witten {\it Two-dimensional gravity and intersection theory on moduli space},
Surveys in differential geometry (Cambridge, MA, 1990), 243--310, Lehigh Univ., Bethlehem, PA, 1991.


\bibitem{Z1} D.~Zvonkine,
{\it An algebra of power series arising in the intersection theory
of moduli spaces of curves and in the enumeration of ramified
coverings of the sphere}, math.AG/0403092~(2004)
%
\end{thebibliography}
\end{document}